\documentclass[11pt,english]{article}
\usepackage[T1]{fontenc}
\usepackage[latin9]{luainputenc}
\usepackage[a4paper]{geometry}
\usepackage{color}
\usepackage{array}
\usepackage{float}
\usepackage{multirow}
\usepackage{amsmath}
\usepackage{amssymb}
\usepackage{graphicx}
\usepackage{rotating}

\makeatletter

\providecommand{\tabularnewline}{\\}

\usepackage{draftwatermark}
\SetWatermarkScale{4}
\SetWatermarkColor[gray]{0.96}

\makeatother

\usepackage{babel}
\usepackage[style=numeric]{biblatex}
\begin{document}
\title{A model for calculating volumes of trains as flows given demand and
capacity restriction}
\author{Martin Aronsson (RISE), Tomas Lid\'en, (VTI)}
\maketitle
\begin{abstract}
This report presents a novel approach to calculate volumes of trains
from a set of requirements describing the transportation demands.
The model presented is based on a multi commodity network flow model.
The model comprises three different layers: the demand level in which
the requirements on the transports are presented, the route layer
which describes the routes possible to implement the demands and the
flow layer in which the resulting flows of train volumes are calculated
and distributed over discrete time periods. This structure corresponds
well to other industrial planning processes where different products
or services that are to be accomplished are first decomposed by methods
into different layers of tasks, by e.g. a work breakdown structure
(WBS). Resources are then allocated to the tasks as well as scheduled
in time. The most important feature of the proposed model is that
it is volume based, i.e. volume of trains are handled instead of scheduling
individual trains. Scheduling a timetable is a very complex task,
and by concentrating on the volumes of trains over time rather than
the detailed scheduling of them, the complexity is greatly reduced
which is a large computational advantage. The disadvantage is that
the representation of capacity in terms of volumes of trains with
different properties becomes crucial. The report discusses these advantages
and disadvantages as well as gives two examples of the use of the
model.

\pagebreak{}
\end{abstract}
\tableofcontents{}\pagebreak{}

The work presented in this report has been supported by the Swedish
Transport Administration (Trafikverket) under grant TRV2019/132113.

\section{\label{sec:Introduction}Introduction}

The following report describes how a operations research model, a
multi commodity network flow model (MCNF), can be used for computing
train volumes from train demands, over a railway network consisting
of stations and links using routing information. The model is intended
to be used in a framework where it is used together with another model
for scheduling TCRs, Temporary Capacity Restrictions. The framework
and the model for scheduling TCRs based on traffic availability is
presented in another parallel report \cite{SATT}. Given named routes
together with volumes entering and exiting the railway network the
MCNF model calculates the train volumes together with the capacity
usage on the link. Crucial for the model to work is the how link capacity
is represented in terms of number of trains per time frame, dependent
on the heterogeneous mix of train types, since no actual scheduling
of individual trains is performed. The solution is a flow of train
volumes over the railway network over discrete time periods. 

The proposed model have three connected layers. On top of the basic
MCNF model are two other layers to handle the fact that different
traffic demands can be implemented by different realizations in terms
of different routes. The topmost layer  is the demand layer, where
demands on volumes of traffic are expressed. The basis for this are
origin-destination pairs (possibly with via-locations) together with
the train type and other basic properties together with the volume
of trains departing in each time period. This layer corresponds to
the demands from railway undertakers (RUs). The middle layer corresponds
to the way the transport is being produced, and corresponds to the
term ``method'' in some other industry sectors. A demand of a large
volume of trains can be implemented by several different routes, thereby
having different ways of producing the demand.

This three layer  approach is found in many industry sectors as well
as in most problem solving methodologies. On a high level  a task
is first split into subtask by methods (describing capabilities) or
break down structures. The subtasks can in turn be further detailed
until a level where all activities are atomic and in some sense ``basically
understood''. Resources are allocated to these activities which are
then scheduled to get a plan such that no resources are overbooked.
This results in a valid plan that accomplish the high end goal. Often
the different activities from the first break down structure are contractually
bound when the should be ready.

The levels are often used in different process steps, typically the
strategic plan is rolling every year with a 36 month coverage, a 24
mount tactical plan and a 12-0 operations process.

The report has three main parts. We first introduce the layered approach
and how it fits into the framework of planning temporary capacity
restrictions as well as a tool for determining a supply or offer of
train volumes. Then the complete mathematical model is presented with
some detailed explanations. Lastly two experimental setups are given,
one smaller to demonstrate some of the features of the model and one
larger to give one of the intended uses of the model.

\section{\label{sec:Multicommodity-flow-network}Multi commodity flow network
flow model}

A multi commodity network flow model (MCNF) \cite{ahuja1993network}
is an operations research model. A network flow model is a directed
graph where each arc has a capacity which must be greater than or
equal to the flow passing that arc. It is also possible to have capacity
on the nodes although we have not elaborated on that in this report.
Network flows are an important special case of linear programming
models, since they (often) preserve the integer property, i.e. if
all source and sink flows are integer, then so are the flows on the
arcs. A MCNF model, however, does not preserve the integer property,
unless the flows of commodities are non-overlapping in the capacity
restrictions in which case each commodity flow problem can be modeled
separately as a network flow problem. In the model presented in this
paper this is not a problem since the volumes of trains are spread
out over several time periods, thus it is not individual trains but
volumes of trains that flows over the arcs.

There are one or more sources (nodes) where flow enters the graph,
and one or more sinks (nodes) where flow is leaving the graph. Each
node must fulfill a balance restriction, the sum of incoming flow
to each node must equal the outgoing flow. Special cases are the sources
and sinks and are modeled differently in the literature. We will regard
incoming and outgoing flows as arc coming from a special kind of nodes
whose only purpose is to act as source or sink. 

MCNF is a network flow model where the flows are of different types,
so called \emph{commodities}. In each node the sum of incoming commodity
must equal the sum of outgoing commodity, so the network keeps track
of the commodities along the arcs and nodes. The source and sinks
are also extended to be of different commodity types. 

The basic formulation of a MCNF problem can be stated as:

\[
\begin{array}{l}
\text{min}\sum_{cij}F(x_{ij}^{c})\\
\text{s.t.}\\
\forall ij\,\sum_{c}x_{ij}^{c}\leq m_{ij}\\
\forall ci\,\sum_{j}x_{ij}^{c}-\sum_{j}x_{ji}^{c}=b_{i}^{k}\\
\forall ijc\,x_{ij}^{c}\geq0
\end{array}
\]

where $x_{ij}^{c}$ is the flow of commodity $c$ on arc $\left\langle i,j\right\rangle $,
$m_{ij}$ is the capacity on arc $\left\langle i,j\right\rangle $
and $b_{i}^{c}$ is either a source node for $c$ (positive value)
or a sink node for $c$ (negative value) or a transit node (zero).

The commodities in the MCNF model presented in this paper are train
services characterized by their origin, destination and type, e.g.
commuter train, long distance passenger train, freight trains. The
are based on what is called Train Service Classes, TSC \cite{Aronsson1302809},
see section \ref{sec:Traffic-Service-Classes}. The network is a graph
based on the original railway network (i.e. stations and lines) and
extended to implement flow of train volumes over time periods. 

\begin{figure}[H]
\begin{centering}
\includegraphics[scale=0.4]{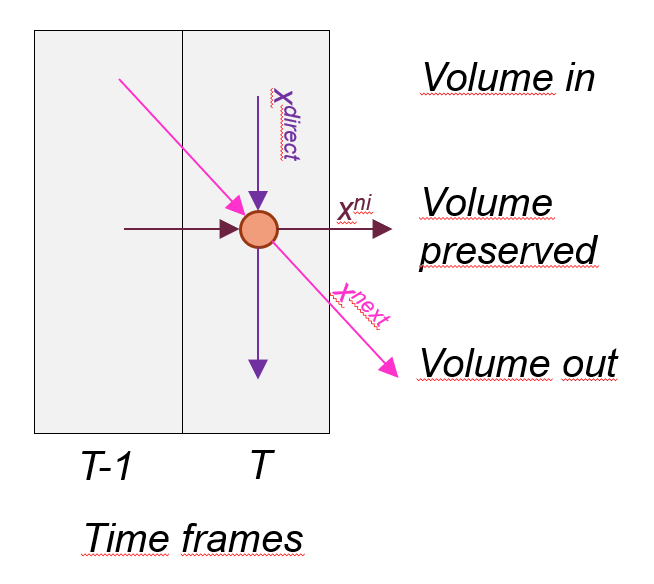}
\par\end{centering}
\caption{\label{fig:Basic-building-block-intro}Basic building block, a node
and its different arcs}
\end{figure}

Let $T$ be a set of time periods $t_{1},...,t_{n}$ then each node
and link of the railway network is transformed into a graph with:
\begin{itemize}
\item $n$ nodes for each station node, one node for each time period. Call
these \emph{timed nodes}.
\item Each resulting timed node is connected with an \emph{inventory arc}
to its consecutive timed node. These hold the volumes that stands
on the station from one time period to the next.
\item Each line between two stations in the railway network is replaced
with $2n-1$ arcs between the same stations, the first $n$ arcs,
called \emph{direct arcs}, are between the stations in each time period,
the next $n-1$ arcs connects consecutive time periods and are called
\emph{next arcs}.
\end{itemize}
in figure \ref{fig:Basic-building-block-intro} this is shown. All
arcs are connected to other nodes, forming a network. 

A small generic example of the network is presented in figure \ref{fig:Basic-flow-model,}
for a railway network consisting of three stations A, B and C and
tracks (directed traffic) A to B and B to C.

\begin{figure}[H]
\begin{centering}
\includegraphics[scale=0.4]{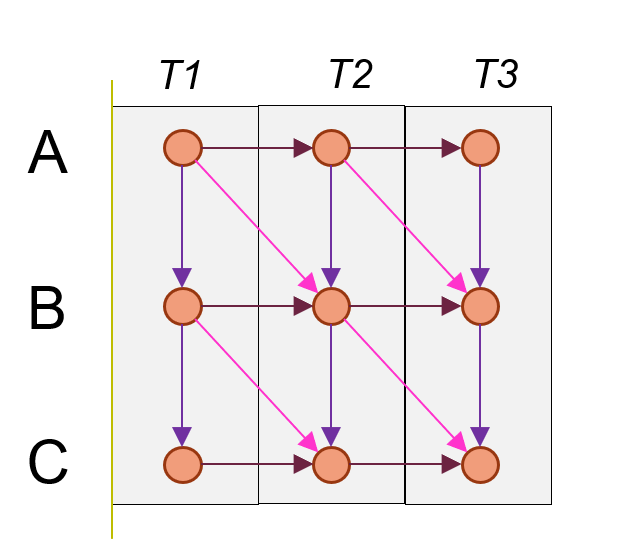}
\par\end{centering}
\caption{\label{fig:Basic-flow-model,}Basic flow model, generic network structure}
\end{figure}

Each line in the geographical network has a total capacity, measured
as a number of units per time period. Certain volumes of train traffic
types consumes more units than other train types, therefore the mix
of different train type volumes are important. The sum of all commodities
on a line in a time period must not exceed the capacity of that line. 

\begin{figure}[H]
\begin{centering}
\includegraphics[scale=0.2]{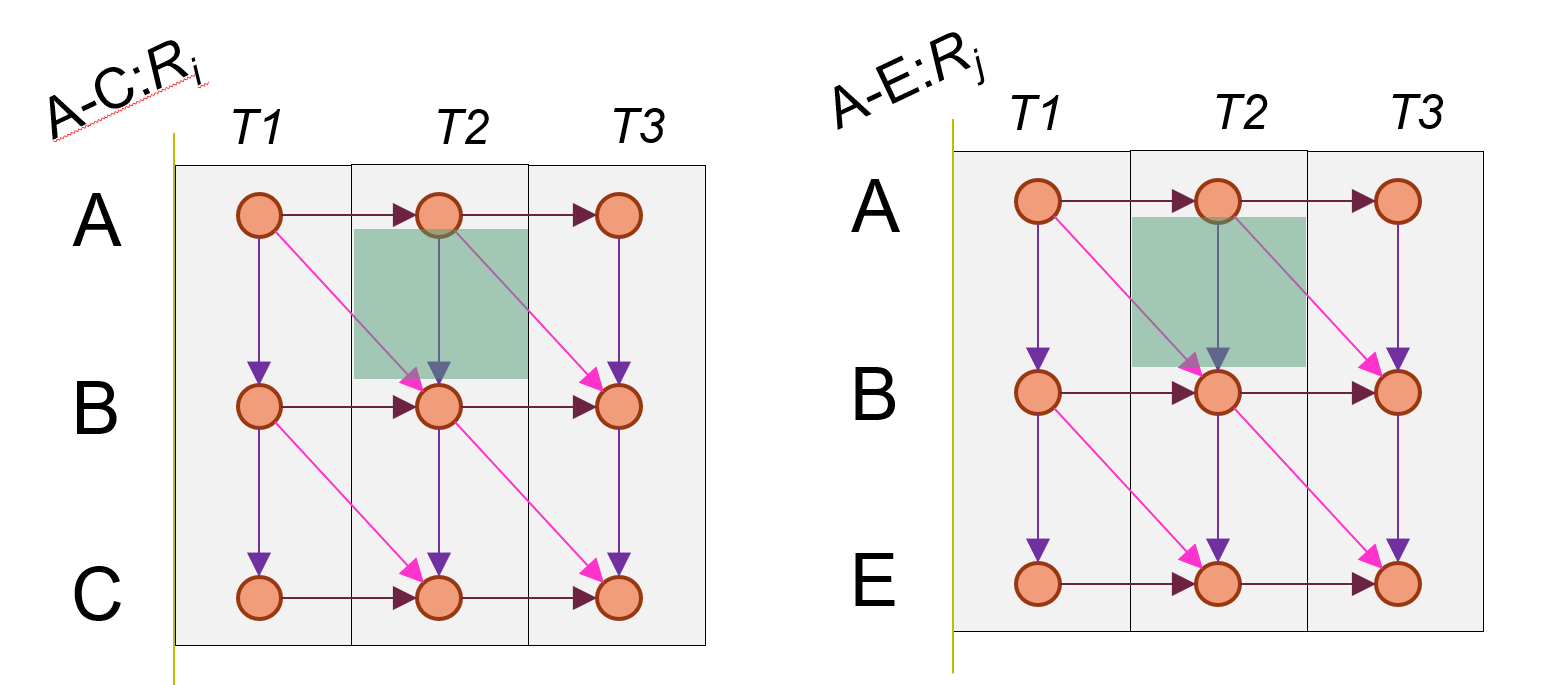}
\par\end{centering}
\caption{\label{fig:Generic-flow-model}Generic flow model with capacity extent}
\end{figure}

Consumed capacity on a line for a time period is calculated as the
sum of all commodities on the direct arcs of that line, half of the
next arcs going into the current time period and half of the next
arcs going out from the current time period. This is illustrated with
the green box in figure \ref{fig:Generic-flow-model}.

This means that the model does not schedule individual trains but
calculates volumes of trains (of different commodities), characterized
by their type, origin and destination, as described in their Train
Service Classes, TSC, in each time period. Different train types can
consume different amount of capacity and the total capacity consumed
must obey the capacity restriction on each link.

Since it is volumes of trains in each time period, in principle no
further details about \emph{when} inside the time period the trains
will depart is known. All that is known is the number of each TSC
(commodity) that will use the line that particular time period. Also
note that the number of units consumed by each TSC is a real number,
indicating the ``spread'' of the actual departure/arrival over time.
The volumes of each TSC that are flowing into the system in each time
period are given (i.e. they are added through the timed source node),
as are the the amount of output flows of the TSC in each timed sink
node. By these ``departure earliest at'' and ``arrives latest at''
it is possible to restrict the volumes of the TSC through the network
so that the volumes actually are flowing according to the train type's
traversal times along the links together with a certain ``spread''
to facilitate meeting and overtakes (The order of the trains are not
scheduled by the flow model, see section \ref{subsec:Traffic-volumes-and-capacity-usage}
which discusses ways of representing capacity as a function of traffic
heterogeneity). Also note that volumes may be split on both the direct
and next arcs making the values on the arcs to be real rather than
integers, making the consumed capacity real.

\section{\label{sec:Traffic-Service-Classes}Traffic Service Classes and Transport
Paths}

\emph{Traffic Service Classes} (TSC) are classes of traffic with equal
properties. The basic properties are the origin, destination, possibly
via-stations, principal performance, possibly axle loads and other
restricting properties, and possibly also time period restrictions
(e.g. not available during certain times of the day). These TSCs are
the commodities, ``products'' or ``package services'' available
to order or apply for. 

A mapping is associated with the TSCs. This mapping describes the
principal ``method'' to implement the packaged service. Each ``method''
consists of a detailed route through the railway net given in the
flow layer  together with different other constraints such as speed
(or, more accurately, the traversal duration for the train type on
the links in the network), axle load and other capabilities of the
TSC. A TSC can also be limited in e.g. when it is available, for example
certain TSCs are only available during rush hours. The idea is that
these TSCs form the ``Product catalog'' of the services available
to implement the demand. Note that several TSCs can implement the
same demand although most commonly there is a best match for each
demand. Also note that each TSC can have several different implementation
methods, i.e. different routes implementing the same service capability
to match the demands.

\begin{figure}[H]
\begin{centering}
\includegraphics[scale=0.3]{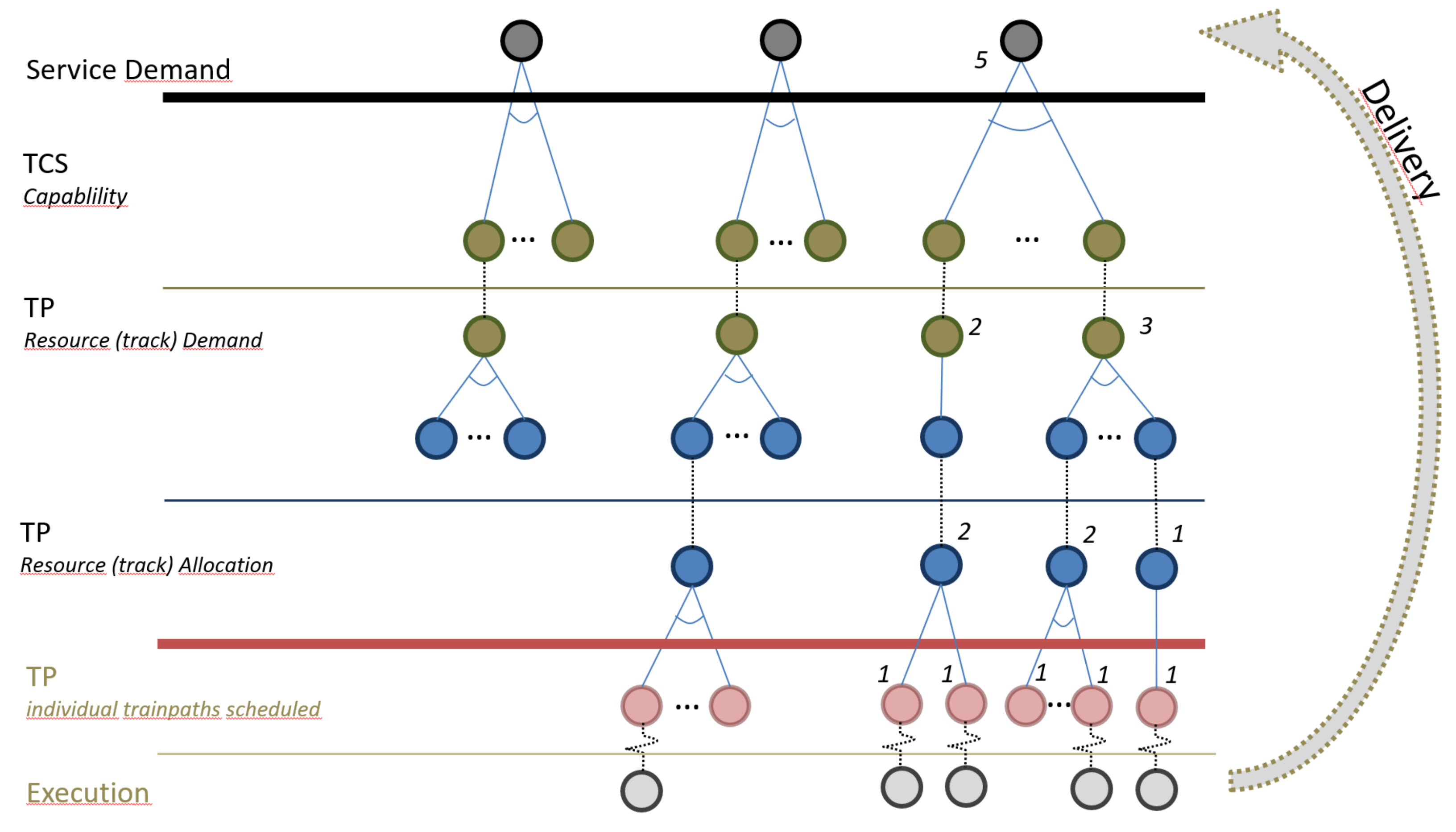}
\par\end{centering}
\caption{\label{fig:Layered-task-structure}Layered task structure}
\end{figure}

So for each given demand, a mapping on a suitable TSC (giving ``qualities''
such as running times etc) is made. This gives a number of possible
routes implementing the given ``offer'' of the TSC. Such an instance
of an TSC is called a \emph{Transport} \emph{path}, TP. The TPs have
all the properties of their TSC together with departure and arrival
timing information as well as a fixed route. Note however that both
choosing a route and finding the way through the flow layer  are solved
together in the proposed model, meaning that the flow layer  affects
the choice of routes and vice versa. In figure \ref{fig:Layered-task-structure}
the overall picture of the different layers and the break down of
a demand into TSCs and TPs is given. 

The model described in this report handles the mapping of demands
down to volumes of trains on links and time periods. The result should
then be used in actual train path scheduling, where individual trains
are scheduled in detail, performed each year to get a valid yearly
capacity allocation. When the detailed capacity allocation plan is
executed, the railway system delivers train transports, which should
be in accordance with the demands original put into the system regarding
requirements, quality and volumes.

\section{\label{sec:The-different-levels}The different layers of service
description}

\begin{figure}[H]
\begin{centering}
\includegraphics[scale=0.4]{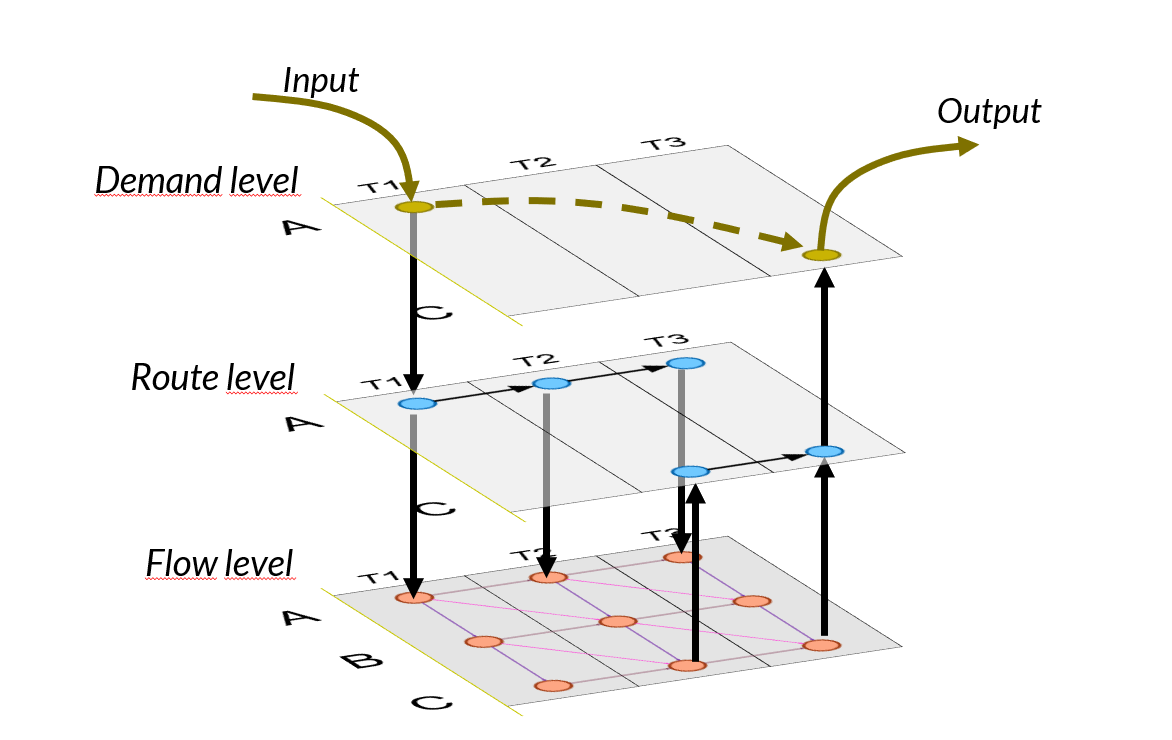}
\par\end{centering}
\caption{\label{fig:Three-layered-task}Three layered task structure used in
the model}
\end{figure}

Figure \ref{fig:Layered-task-structure} above describe the principle
break down of demands to train paths. We will now detail the different
layers using three dimension figures, complemented by graphs. The
three dimensions are the railway network, the layers of abstraction
and time discretized into periods of equal length. At the demand layer
 are the sources and sinks of the flow problem, see figure \ref{fig:Three-layered-task}.

All the layers are formulated together in the same model and computation
is performed on the combined model. Note that the time periods are
the same for all layers. It should be possible to have a more fine-grained
time scale in the flow model compared to the demand layer, but we
will not investigate that further in this report. 

The different layers are described separately since the correspond
to different perspectives of the problem formulation. The vertical
arcs (arrows) between the different abstraction layers are not arcs
in the flow model and should be understood as communication between
the layers and in the concrete model often the same node. They are
separated here to facilitate the understanding the sub-problems making
up the complete model.

\subsection{\label{subsec:Demand-level}Demand layer}

\begin{figure}[H]
\begin{centering}
\includegraphics[scale=0.2]{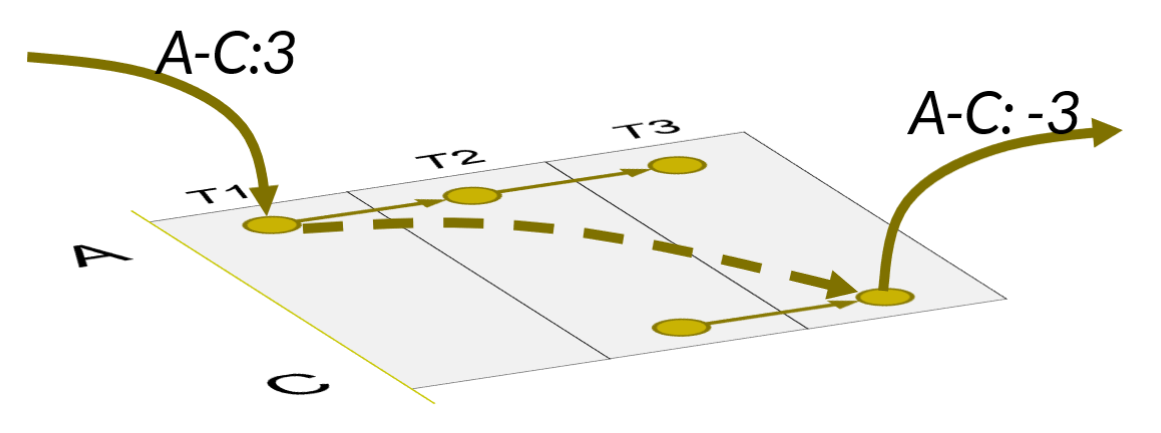}
\par\end{centering}
\caption{\label{fig:Demand-layer}Demand layer}
\end{figure}

The demand layer  (DL) is where the traffic demand is given. The volumes
entering at DL is given as a demand in terms of the number of trains
of each TSC. In figure \ref{fig:Demand-layer} a demand of 3 trains
of demand A-C between stations A and C is given, with the timing constraints
that all three trains should depart in time period $T_{1}$ and be
at station $C$ at time period $T_{3}$.

Thus the named demand $A-C$ is described as a origin-destination
pair (possibly with via-stations) together with the volumes departing
(arriving) in each time period in the figure. This can be stated,
given in semi-mathematical notation, as ``input(A-C) in $T_{1}$:3
units'' and ``output(A-C) in $T_{3}$:3 units'' shown by the arrows
going into and out of the graph. The dashed line shows what is sought
in the lower layers of the model for the demand to be realized, i.e.
the sought implementation of the demand. Note that demand is given
as a volume in each time period, not as individual train path requests.
This means that the individual departures are equally probable in
the time period, i.e. it is not possible to state a specific time
when the train is to depart. A core feature of the model is that it
is \emph{volumes} of traffic that is sought, not individual trains.
However, by making the time periods small enough it is possible to
get more detailed departure time domain. It should also be possible
to give the departure volumes a certain probability ``shape'' over
the time period, although this is not elaborated further in this report.

The demand layer  is also where volumes can be delayed from the applied
demanded time of departure, which are the arrows in the right picture
going from one time period to the next at the origin (analogously
defined at the destination, if there are requirements on not arriving
too early). A cost is associated with delaying the departure. Note
that when some of the volume is delayed from its origin this volume
is not standing at the station (A) waiting to depart. Rather this
is breaking the demand requirement and should be punished by the objective
function (whereas if the volume would be standing at the origin station
in the flow layer  it is consuming capacity at the station waiting
to depart). With these route layer  delays it is possible to capture
that the offer to the demand has a discrepancy. It is different to
offer another departure time of the train compared to taking the demanded
departure time and then letting the train wait to depart at the station.

\subsection{\label{subsec:Route-level}Route layer}

Routes are implementations of demands. There could be several routes
matching a demand, and part of the problem solving is to choose which
routes that implement the demands. Different routes implement the
demands differently, with different drawbacks, and the objective function
contains the costs for choosing a particular route (as well as costs
for delays etc).

\begin{figure}[H]
\begin{centering}
\includegraphics[scale=0.2]{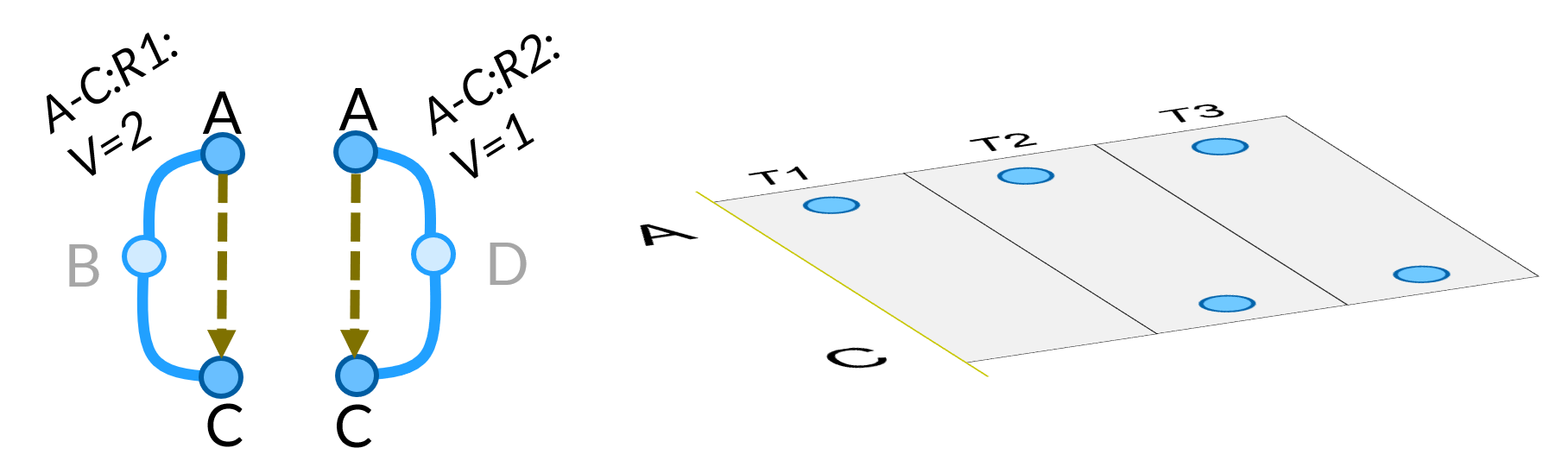}
\par\end{centering}
\caption{\label{fig:Route-layer}Route layer}
\end{figure}

The route layer  takes care of choosing the route given the demand
from the previous layer, through the (geographical) network. The mapping
is performed by the TSCs which as part of its description hasa number
of different implementations of the demand by different routes. Each
route has a transport duration associated with it together with other
properties such as axle load, time of day it is applicable etc. The
TSCs' mapping of demands are made on the basis of each demand's requirements,
e.g. the origin-destination pair is one of the important ones. In
figure \ref{fig:Route-layer} there are two TSCs having the capability
to match the demand of volumes starting in A and finishing in C. The
two possible routes to implement the demand A-C are R1 going A-B-C
and R2 going A-D-C, named A-C:R1 and A-C:R2 in figure \ref{fig:Route-layer}.

Routes correspond to method descriptions in other industry sectors
and are an important tool in planning. Methods are named ways of performing
a task, having different properties such as e.g. execution time and
reource consumption often as a function of amounts of e.g. materials.
By choosing methods for one, several or all tasks of a project different
make spans results together with different total costs.

\subsection{\label{subsec:Flow-level}Flow layer}

Most of the model is formulated around the multi-commodity network
flow model (MCNF), and it is on this layer  the source and sink of
the flow are connected. Time has the same discretization into a number
of time periods as previous layers. Each TSC and their corresponding
route description corresponds to a commodity in the MCNF problem formulation
called a named route (i.e. A-C:R1 and A-C:R2 in figure \ref{fig:Route-layer}).
The flow graph is built up as follows.

In each time period every station node is represented as a node. From
every station node there are three different types of arcs that can
go into the node and three different types of arcs that can go out
of the node. The three types of arcs are:
\begin{itemize}
\item \emph{node inventory arcs}, one for each possible named route, representing
volumes of each named route standing at the station between two consecutive
time periods,
\item \emph{direct arcs}, one for each named route, going from one node
to another node along a railway network link. Direct arcs hold volumes
of each named route traversing the link in the same time period
\item \emph{next arcs}, one for each named route, going from one node to
another along a railway network link Next arcs take volumes of each
named route traversing the link from one time period to the next time
period.
\end{itemize}
There may be zero, one or more arcs of each type going into or out
of the node. This is visualized in figure \ref{fig:Basic-building-block,}. 

\begin{figure}[H]
\begin{centering}
\includegraphics[scale=0.4]{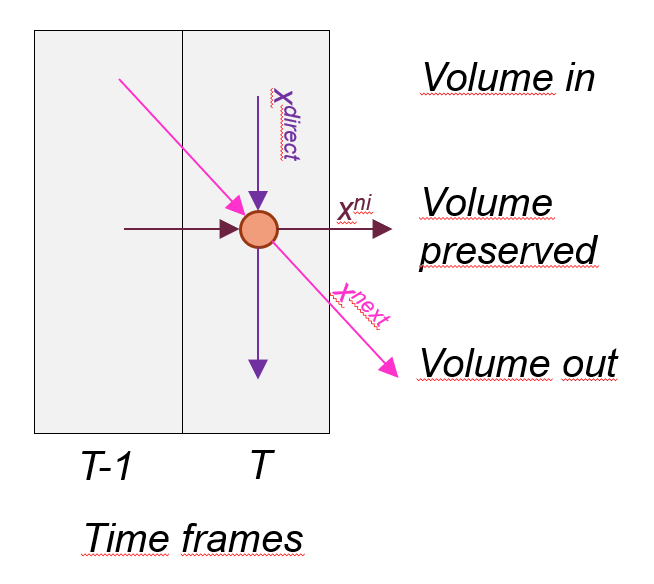}
\par\end{centering}
\caption{\label{fig:Basic-building-block,}Basic building block constiting
of a node and all of it's different arc types going into and out of
the node}
\end{figure}

Note that all these arcs take real number values, i.e. they do not
have to be integer, since the demand is given as a equally distributed
volume over the time period. As this volume is moving through the
net, fractional parts of that volume is e.g. taking next arcs to the
next time frame. 

\begin{figure}[H]
\begin{centering}
\includegraphics[scale=0.25]{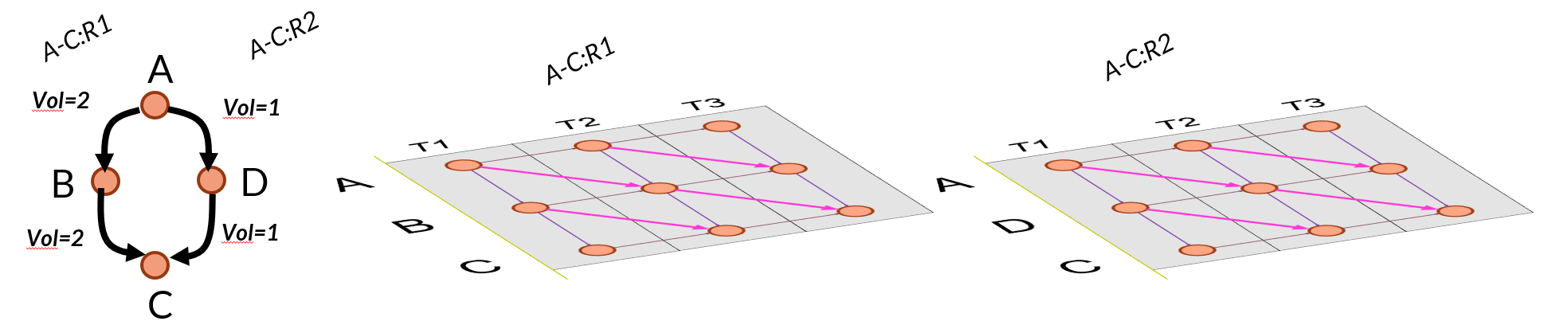}
\par\end{centering}
\caption{\label{fig:Two-route's-implementation}Two route's implementation
in the flow layer}
\end{figure}

For each named route we have a network of timed station nodes and
track links as showed in figure \ref{fig:Two-route's-implementation},
where the sub-graphs of routes R1 going A-B-C and route R2 going A-D-C
are both given. 

Many named routes can use the same network link, so for each time
period $T_{n}$ and network link all arcs for a each named route are
added together to get the capacity usage of that network link and
time period $T_{n}$. Within a time period at most a maximum number
of trains can be accommodated on that link, i.e. the capacity usage
must be withing the capacity limit for that network link. To calculate
the capacity usage in a link. all direct arcs for each named route
is summarized, together with half of the sum of the incoming next
arcs and half of the sum of the outgoing next arcs, for all named
routes. This is shown in figure \ref{fig:Basic-mode,-capacity} for
two named routes A-C:$R_{i}$ and A-E:$R_{j}$. 

\begin{figure}[H]
\begin{centering}
\includegraphics[scale=0.2]{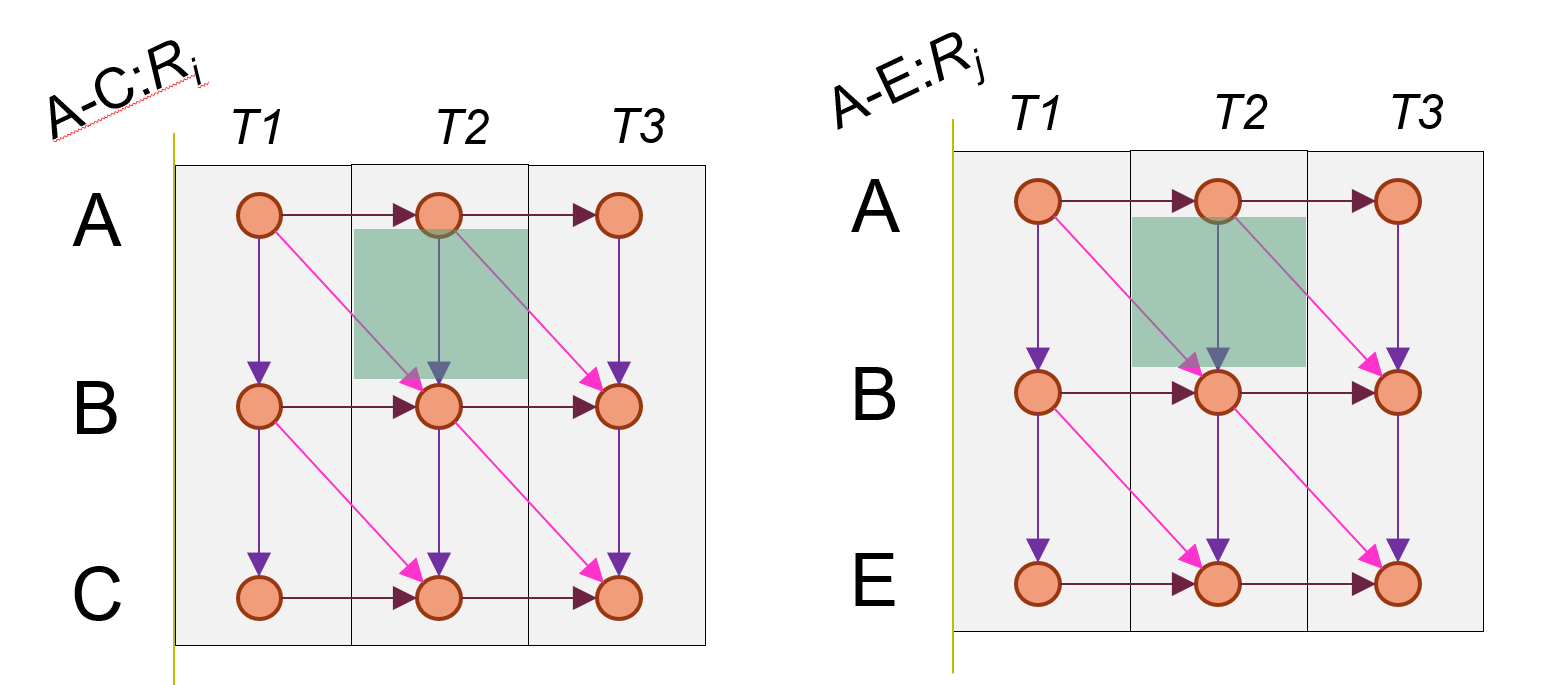}
\par\end{centering}
\caption{\label{fig:Basic-mode,-capacity}Basic mode, capacity extent}
\end{figure}

Taking half of the incoming and half of the outgoing arcs' values
is an approximation. One motivation for that is that the next arc's
volumes, given that the volume is evenly distributed over the time
period, has half of its capacity consumption inside time frame $T_{n}$
and half in the next time frame $T_{n+1}$.

Note that the values on the arcs are fractional. This is natural since
we are ``scheduling'' volumes of trains that depart somewhere inside
the time period $T$, not train individuals, and it is parts of volumes
that passes on to the next time period, i.e. parts of the evenly distributed
volume from time period $T$, more about this in section \ref{subsec:Traffic-volumes-and-capacity-usage}.
Therefore it is not strange that there can be a fractional number
of trains taking an arc. It is the capacity consumption that the model
measures, not the scheduling of individual trains.

The flow part adds all volumes together on the timed station and track
nodes and arcs, where each of them has a capacity limit. The split
of volumes in named routes is made at the route layer, but the ``flow
over time'' is determined at the flow layer. So each volume departing
at a certain station node the flow layer  determines the flow to the
arrival station.

\subsection{\label{subsec:Traffic-volumes-and-capacity-usage}Traffic volumes
and capacity usage on line segments in the flow layer}

\begin{figure}[H]
\begin{centering}
\includegraphics[scale=0.5]{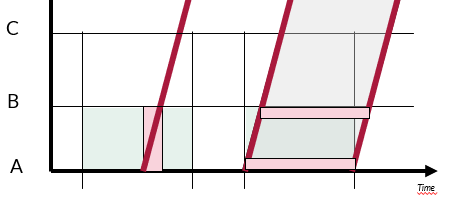}
\par\end{centering}
\caption{\label{fig:Train-paths-and}Train paths and volumes}
\end{figure}

In this flow model, traffic is modeled as a volume that arrives/departs
inside a time period. This means that there are no events in the model
when a train arrives or departs, but volumes that occur within a time
period. The actual departure time is not an element of the flow model,
but the time period where the train departs, and the probability for
when, inside the time period, that the train will depart is equally
distributed. The time (capacity) the trains will occupy the line
segment is transformed to a ``capacity volume'' that the trains
will allocate of that line segment's total capacity in that time frame.
This is illustrated in figure \ref{fig:Train-paths-and} where the
pink bars have the same area as well as the light green areas (the
capacity of section A-B). 

Volumes are moving along its path (A to B), and when it does so, the
volumes are forced to gradually move into the next time period, in
line with the speed of the original train (we omit the possibility
for volumes to take the inventory arcs in this section). If there
are four trains of a specific TSC departing within a time period,
the volume added at the source in the flow model corresponding to
the origin (station) and time period will be 4 inside this particular
time period. If the time period is one hour and the travel time of
the TSC over the line segment A-B is 15 minutes, then a volume of
4 trains leaving A that would have passed into the next time period
when arriving at B is 1

\begin{figure}[H]
\begin{centering}
\includegraphics[scale=0.5]{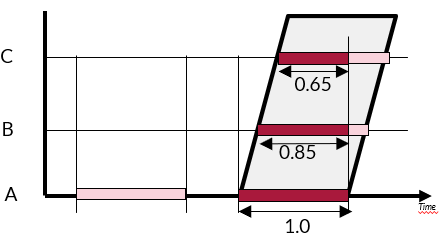}
\par\end{centering}
\caption{\label{fig:Volumes-moving-along}Volumes moving along route}
\end{figure}

Volumes may therefore be real numbers as they travel along their paths.
A train that departs from its origin A will have a volume of 1, equally
distributed over the time period. As this equally distributed volume
moves along the line segment it will consume time, and therefore a
fraction of it will be present in the next time period. This is illustrated
in figure \ref{fig:Volumes-moving-along} where 15\% of the volume
has moved into the next time period for line segment A-B, and 20 \%
when moving from B to C.

\begin{figure}[H]
\begin{centering}
\includegraphics[scale=0.5]{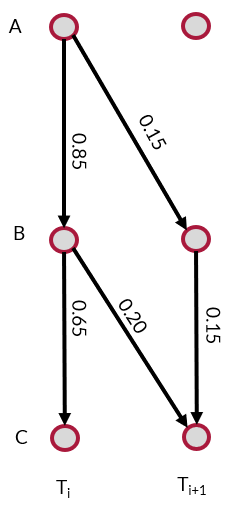}
\par\end{centering}
\caption{\label{fig:Flows-on-arcs}Flows on arcs in mode}
\end{figure}

In the flow model (only the part shown in the previous figures and
omitting the node inventory arcs) the nodes and arcs will be as in
figure \ref{fig:Flows-on-arcs}.

\begin{figure}[H]
\begin{centering}
\includegraphics[scale=0.5]{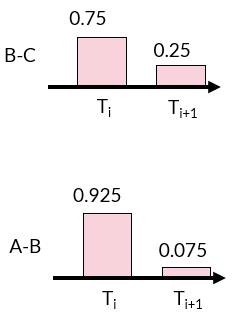}
\par\end{centering}
\caption{\label{fig:Volumes-on-arcs}Volumes on arcs over time}
\end{figure}

The capacity used by a volume of a TSC over a line segment is the
sum of the direct arcs plus half of the next arcs going to the next
time period as well as half of the volumes coming into the current
time period. In the example there are no incoming volumes to the time
period $T_{i}$ so the capacity used in time period $T_{i}$ is $0.85+0.15/2=0.925$
capacity units. For the line segment BC during time period $T_{i+1}$
it is $0.2/2+0.15=0.25$. Capacity usage bar charts are given in figure
\ref{fig:Volumes-on-arcs}. 

Corresponding bar charts can given for the nodes (stations). These
are important since these will be used in the flow model to calculate
the maximum volumes that can possibly have arrived to a stations each
time period. These sums are used to make the flows take the ``next''
arcs not only the direct arcs. These sums are in this paper called
'maximum aggregated sums' and are calculated for each node, time period
and TSC. So for example, in period $T_{i}$ at node C the maximum
volume of the example TSC above can at most be $0.65$, the other
0.35 must have moved into time period $T_{i+1}$. 

\begin{figure}
\begin{centering}
\includegraphics[scale=0.25]{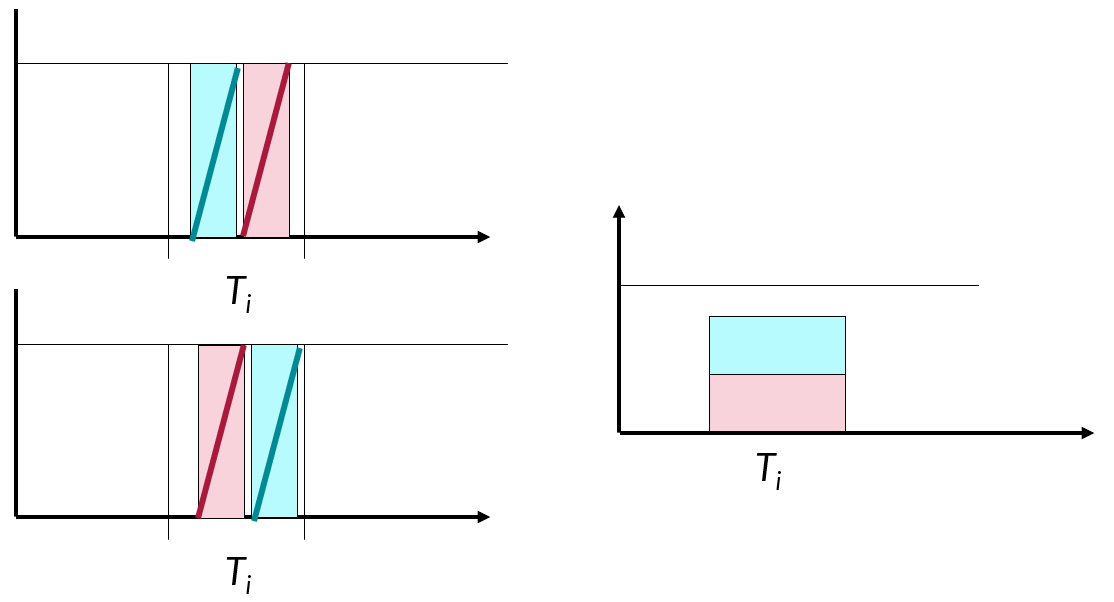}
\par\end{centering}
\caption{\label{fig:Within-time-period}Within time period order of trains
have no meaning}
\end{figure}

If other TSCs also are using the line segment, the capacity used is
the sum of those TSCs. The bars in the bar chart has to be added and
should be below the capacity limit of the line segment and time period.
Since the time is discretizised into periods, the order of the instances
of the different TSCs added is not known, only that the capacity consumption
will be the cumulative sum of the two, see figure \ref{fig:Within-time-period}.

\begin{figure}[H]
\begin{centering}
\includegraphics[scale=0.25]{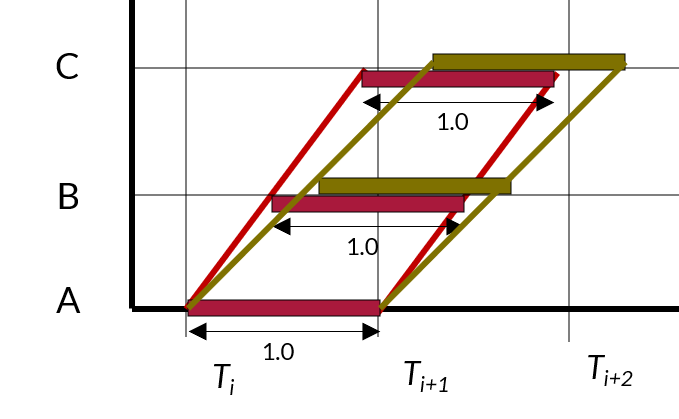}
\par\end{centering}
\caption{\label{fig:Maximum-and-minimum}Maximum and minimum aggregated sums}
\end{figure}

An analogous measure at each node and time frame can also be introduced,
the minimum volumes that must have reached the node. This put emphasis
on that volumes are actually traveling through the net at a particular
pace. The 'minimum aggregated sum' says that the efficiency (volumes
arrived at a station at a certain time period) must be larger than
a threshold, at each node. This measure is of more importance at the
sinks of the flow graph. the importance of this measure is that volumes
cannot be delayed too much early in the flow so that the volumes cannot
reach the sink in time for its latest arrival. This is shown in \ref{fig:Maximum-and-minimum}.

\begin{figure}[H]
\begin{centering}
\includegraphics[scale=0.25]{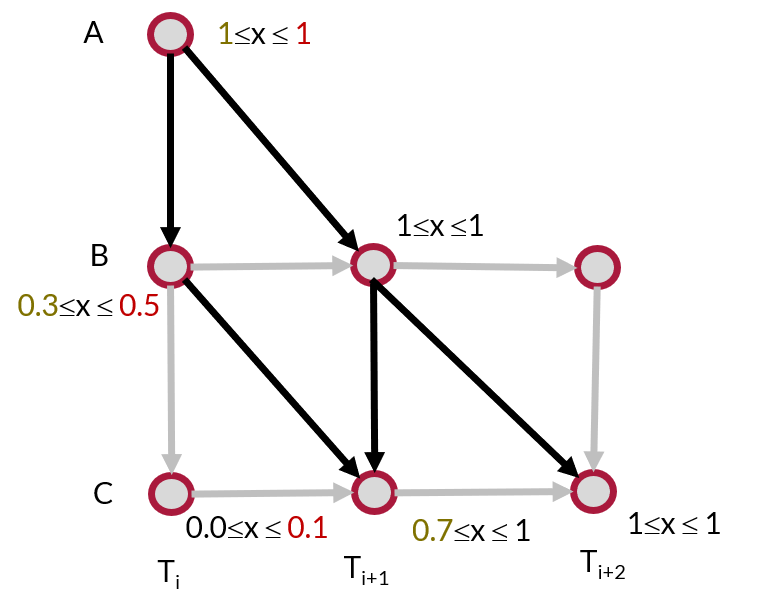}
\par\end{centering}
\caption{\label{fig:Minimum-and-maximum}Minimum and maximum aggregated sums
in nodes}
\end{figure}

Green bars are the minimum pace the TSC must have, which has less
speed than the best performance shown in red. The corresponding flow
graph with the aggregated sums at the nodes is shown in figure \ref{fig:Minimum-and-maximum}
where the minimum and maximum numbers are shown in the same color
as in figure \ref{fig:Maximum-and-minimum}.

When we use the term 'aggregated sum' we will refer to the 'maximum
aggregated sum', since this is crucial for having the model work the
correct way so that volumes are ``pushed'' forward in time. The
maximum aggregated sums make sure the volumes does not break the maximum
performance of their corresponding trains. The minimum aggregated
sum is useful for stating efficiency constraints for the TSCs, and
to uphold the latest arrival time at the sink.

\subsection{\label{subsec:Capacity}Capacity}

Network capacity is handled on the flow layer. The basic approach
in this paper is that capacity is measured as the sum of all instances
of TSC's in each time frame at each link (as described in section
\ref{subsec:Flow-level}).

\begin{figure}[H]
\begin{centering}
\includegraphics[scale=0.5]{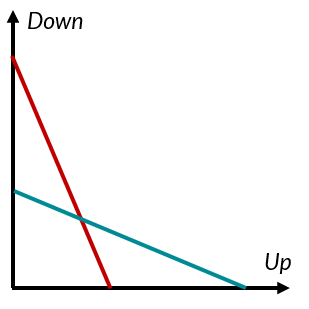}
\par\end{centering}
\caption{\label{fig:Principal-model-for}Principal model for setup times for
single track lines, the y axis is capacity in down direction while
x axis represents the capacity for the up direction }
\end{figure}

Double track lines are represented with two directed network arcs
in opposite direction. Single track lines are represented by two arcs
too, one in each direction. However, since they correspond to one
physical track in reality, there is a function to merge the two volumes
on the two network arcs in such a way that the two volumes can be
realized on the same physical track despite the fact that the move
in opposite direction. Change of direction is handled as a setup time
between TSCs in opposite directions. But as the model does not schedule
the TSCs inside a time frame, a function must be introduced that returns
a probable sum of the setup times that is likely to be needed. This
function is preferable sensitive to the mix of how many TSCs there
are in the two directions. The more even the two directions are and
with higher volumes, the more setup time is need. This gives us a
concave function, shown in figure \ref{fig:Principal-model-for}.

If there are only red TSCs moving in the 'Down' direction then the
nominal capacity could be utilized (as is the case for double track
lines). But with more blue TSCs moving in the 'Up' direction, the
number of red TSCs are decreasing faster up to the point where there
are as many blue and red TSCs (which is when the line has its worst
performance). To model this an integer variable has to be introduced
which represents the point where the blue and red lines meet in figure
\ref{fig:Principal-model-for}. The constraint is generally formulated
as four equations, where $x$, $y$ are the volumes in each direction,
$x_{c}$, $y_{c}$ are the capacity limit in each direction and $C_{max}$
is the capacity limit (here equal in both directions). $K\geq1$ is
meeting setup coefficient, if $K$ is chosen high, the setup time
for changing direction is high, if it is low the setup time is less.
When $K=1$ there is no setup time.

\[
\begin{array}{l}
x_{c}=C_{max}-K\times y\\
0\leq x\leq x_{c}\\
y_{c}=C_{max}-K\times x\\
0\leq y\leq y_{c}
\end{array}
\]

For double track lines, no setup time for changing direction is used.
As links in the model have direction, a single track line is represented
as two links. This means that there are two logical links (one in
each direction) using the same physical link. This means that the
capacity calculations in the model, for single track lines, 

There is a nominal capacity on each link, based on the most common
capacity consumption for the common TCSs. All TCSs then have their
own capacity limit (a variable). The sum of all TCSs own capacity
limits together with all setup times and together with capacity consumption
for TCRs shall not exceed the nominal capacity for that link. For
TSCs with deviating duration (faster or slower) compared to the nominal
TSC's duration an extra capacity consumption factor is applied to
the scheduled number of that TCS.

\section{\label{sec:Mathematical-model}Mathematical model}

This section introduces the mathematical presentation of the model.
It is important to note the distinction between the geographical nodes
and links, called station nodes $N$ and track links $L$ representing
the geographical network, and the problem graphical formulation where
the nodes are $N\times T$ i.e. all stations nodes at all time periods
and $L\times T$ which are the arcs also for all time periods, i.e.
direct and next arcs. Between consecutive time periods are also node
inventory arcs in $N\times T_{i}\times T_{i+1}$ drawn. 

\subsection{\label{subsec:Notation}Notation}

We use Greek letters for functions (such as $\sigma)$, hollow capital
letters for sets (such as $\mathbb{D}$), bold capital letters for
matrices and multi-dimensional arrays (such as $\boldsymbol{C}$).
Lowercase letters are used for elements of sets and matrices (such
as $l$) as well as for indexes.

Functions will be use for mappings, and when mapped onto tuples with
more than one element will often use a selection function $\pi^{i}$
to reference the $i$:th element of the tuple. As an example take
the map from demand names onto a tuple of its properties:

\[
\tau:\mathbb{D}\mapsto\{\left\langle n_{o},n_{d},h\right\rangle |n_{o}\in\mathbb{N},n_{d}\in\mathbb{N},h\in\mathbb{H}\}
\]

If $\tau(d)=\left\langle n_{1},n_{2},h\right\rangle $, some $n_{1}$
and $n_{2}$ . Then, writing

\[
\forall d\in\mathbb{D},\exists h\in\mathbb{H}:\pi^{3}(\tau(d))=h
\]

is the same as

\[
\forall d\in\mathbb{D},\exists h\in\mathbb{H},\exists n_{1}\in\mathbb{N},\exists n_{2}\in\mathbb{N}:(\tau(d)=\left\langle n_{1},n_{2},h\right\rangle \wedge\pi^{3}(\left\langle n_{1},n_{2},h\right\rangle )=h))
\]

\subsubsection{\label{subsec:Enumerations}Enumerations}

We will use the following basic sets of atomic elements. The basic
types are enumerated, i.e. each element of the set is assigned an
unique integer (a unique index), and the largest integer in the enumeration
is the same as the cardinality of the set. Enumerated sets will have
a hollow capital letter, such as $\mathbb{R}$. With slight abuse
of notation we will use the notation $\mathbb{R}_{i}$ when referencing
the $i$:th value of the enumeration of the set.. We use lower case
letters to reference elements in the set, so $\forall r\in\mathbb{R}\exists i:r=\mathbb{R}_{i}$
states that for all elements $r$ in the set $R$ there exists an
index $i$ such the the $i$:th element of $\mathbb{R}$ is $r$.

\paragraph{\label{par:Basic-types}Basic types}
\begin{description}
\item [{$\mathbb{T}$}] A set of time periods, numbered from 1 and up to
the largest investigated time period $t_{max}$. We often use the
letter $t$ to denote a time period. 
\item [{\textmd{$\widehat{\mathbb{T}}$}}] $\mathbb{T}$ extended with
the time period 0, i.e. initial values in the model (e.g. initial
volumes on nodes).
\item [{$\mathbb{N}$}] The enumerated set of station nodes $n\in\mathbb{N}$
in the investigated network. We will use the notation $n_{i}$ to
reference the $i$:th element for the enumerated set $\mathbb{N}$.
\item [{$\mathbb{L}$}] The enumerated set of track links $l\in L$ in
the network. 
\item [{$\varPhi$}] Mapping of link names to its properties, $\varPhi:\mathbb{L}\mapsto\{\left\langle n_{i},n_{j}\right\rangle |(n_{i},n_{j})\in\mathbb{N}^{2}\text{ and }n_{i}\neq n_{j}\}$A
link goes from node $n_{i}$to $n_{j}$ and has a direction. Moreover,
$n_{i}$ is called the \emph{tail} of the link $l$ , $\pi^{1}(\varPhi(l))=n_{i}$
and $n_{j}$ is called the \emph{head} of the link $l$, $\pi^{2}(\varPhi(l))=n_{j}$.
A single track line has two links, one in each direction. 
\item [{$\mathbb{H}$}] The enumerated set of train types. We will use
$\{reg,ic,gt\}$ standing for regional trains, intercity trains and
freight trains.
\item [{$\mathbb{D}$}] The enumerated set of demands (demand names) for
transports from an origin $n_{i}$ to a destination $n_{j}$. As for
links, a demand has a direction and goes from one node (the origin)
to another node (the destination) in the network.
\item [{$\mathbb{R}$}] The enumerated set of route names for all paths
through the net. The actual routing through the network is given by
the matrix $\boldsymbol{R}^{L}$ described below.
\end{description}

\subsubsection{\label{subsec:Arrays-(lists),-matrixes}Arrays (lists), matrices
and other data structures of ground values}

The basic data types are ground (non-variable) and given as part of
the problem formulation.

The following data structures models the production (flow) layer  of
section \ref{subsec:Flow-level}.
\begin{description}
\item [{$\sigma$}] coupled links, $\sigma:\mathbb{L}\mapsto\mathbb{L}$.\\
Function mapping every $l\in\mathbb{L}$ to an element $l\in\mathbb{L}$.
$\sigma$ is used to distinguish between single track lines and double
track lines. For double track lines, $\sigma(l_{i})=l_{i}$ i.e. $l_{i}$
is mapped to itself if it is part of a double track line, while for
single track lines $\sigma(l_{i})=l_{j}\text{ where \ensuremath{\pi^{1}(\varPhi(l_{i}))=\pi^{2}(\varPhi(l_{j}))\wedge\pi^{1}(\varPhi(l_{j}))=\pi^{2}(\varPhi(l_{i}))}}$,
i.e. a single track line in $\mathbb{L}$ is by $\sigma$ mapped onto
its reverse link (using the same physical track).
\item [{$\boldsymbol{C}_{lt}^{nom}$}] Matrix of nominal capacities (real
number) ranges over $\mathbb{L}\times\mathbb{T}$.\\
The nominal capacity is the largest volume of trains that can pass
over a time period. The nominal capacity can be restricted by possessions,
heterogeneous traffic etc. Nominal capacity is given based on the
typical train type, and is decreased by factors for heterogeneous
traffic (i,e, train types other than the chosen typical train type).
\item [{$\boldsymbol{D}_{lh}$}] Matrix of real values ranging over $\mathbb{L}\times\mathbb{H}$.\\
The duration a train of type $h\in\mathbb{H}$ takes to traverse the
link $l\in\mathbb{L}$.
\end{description}
The following data structures concerns the demand for traffic given
at the demand layer of section \ref{subsec:Demand-level}. The index
set $\mathbb{D}$ contains the ``names'' of the offered services
and $\tau$ maps the names onto its properties. The matrix $\boldsymbol{D}^{m}$
relates demand names with route names, this relation is a one to many
relation. The matrix $\boldsymbol{D}^{V}$ gives the demanded traffic
volumes over time.
\begin{description}
\item [{$\tau$}] Mapping from demand names to its properties origin, destination
and train type, $\tau:\mathbb{D}\mapsto\{\left\langle n_{o},n_{d},h\right\rangle |n_{o}\in\mathbb{N},n_{d}\in\mathbb{N},h\in\mathbb{H}\}$
\item [{$\boldsymbol{D}_{dr}^{m}$}] Matrix of binary values ranges over
$\mathbb{D}\times\mathbb{R}$.\\
When $\boldsymbol{D}_{dr}^{m}=1$ then route $r\in\mathbb{R}$ implements
(realizes) the demand $d\in\mathbb{D}$. There may be many alternative
$r_{j},\,...,\,r_{k}$ that implements a demand $d\in\mathbb{D}$.
\item [{$\boldsymbol{D}_{dt}^{V}$}] Matrix of integers ranging over $\mathbb{D}\times\mathbb{T}$.\\
$\boldsymbol{D}_{dt}^{V}$ contains the demanded volumes of each member
$d\in\mathbb{D}$ in each time period $t\in\mathbb{T}$.
\item [{$\nu$}] $\nu:\mathbb{D}\mapsto\mathfrak{R}$ Function that maps
the sum of demands over time to a real number, defined as: $\forall d\in\mathbb{D}:\nu(d)=\sum_{t\in\mathbb{T}}\boldsymbol{D}_{dt}^{V}$
\end{description}
The following data structures concerns the possible routes $r\in\mathbb{R}$
that implement one or several demands from $\mathbb{D}$ and corresponds
to the routing level of section \ref{subsec:Route-level}. The content
of this routing layer corresponds to the TSCs, Traffic Service Classes,
described in section \ref{sec:Traffic-Service-Classes}. The routes
are also the commodities of the multi-commodity flow model.
\begin{description}
\item [{$\psi$}] Mapping from route names to its properties origin, destination
and train type, $\psi:\mathbb{R}\mapsto\{\left\langle n_{o},n_{d},h\right\rangle |n_{o}\in\mathbb{N},n_{d}\in\mathbb{N},h\in\mathbb{H}\}$
\item [{$\boldsymbol{R}_{rl}^{L}$}] Matrix of binary values ranging over
$\mathbb{R}\times\mathbb{L}$.\\
If route $r$ uses link $l$ then $\boldsymbol{R}_{rl}^{L}=1$ else
$\boldsymbol{R}_{rl}^{L}=0$. This means that a route can only use
a link once (i.e. there may not be cycles in the route). 
\end{description}
Two links $l_{i},l_{j}\in\mathbb{L}$ in the path of a route $r\in\mathbb{R}$
are said to be ordered, $l_{i}\prec_{r}l_{j}$, if $\boldsymbol{R}_{rl_{i}}^{L}=\boldsymbol{R}_{rl_{j}}^{L}=1\wedge\pi^{2}(\varPhi(l_{i}))=\pi^{1}(\varPhi(l_{j}))$.
We say that $l_{i}$ \emph{precedes} $l_{j}$ in the path of route
$r$, and that $l_{j}$ is the \emph{following link} of $l_{i}.$Note
that this implies that a link can only be passed once in a route. 
\begin{description}
\item [{$\boldsymbol{R}_{rn}^{N}$}] Matrix of binary values ranging over
$\mathbb{R}\times\mathbb{N}$.\\
If route $r\in\mathbb{R}$ uses node $n\in\mathbb{N}$ then $\boldsymbol{R}_{rn}^{N}=1$
else $\boldsymbol{R}_{rn}^{N}=0$. 
\item [{$\mathbb{I}(o)$}] Define the indexed set $\mathbb{I}(o)$ indexed
by the demand $o\in\mathbb{D}$ as\\
$\forall o\in\mathbb{D}:\tau(o)=\left\langle n_{1},n_{2},h\right\rangle \wedge\mathbb{I}(o)=\{r|r\in\mathbb{R}\wedge\psi(r)=\left\langle n_{1},n_{2},h\right\rangle \}$
\\
i.e. $\mathbb{I}(o)$ consist of all routes $r$ that implements $o$.
\\
If $\mathbb{I}(o)=\emptyset$ the demand $o$ cannot be performed. 
\item [{$\boldsymbol{D}_{rn}^{aggr}$}] Matrix of real number values ranging
over $\mathbb{R}\times\mathbb{N}$.\\
$\boldsymbol{D}_{rn}^{aggr}$ is the \emph{maximum aggregated sum}
of all transport durations on each link up to the node $n\in\mathbb{N}$
in the route $r\in\mathbb{R}$, i.e. let $\pi^{1}(\varPhi(l))=n$
and $\psi(r)=\left\langle n_{1},n_{2},h\right\rangle $ for some $n_{1},n_{2}$,
then $\boldsymbol{D}_{rn}^{aggr}=\sum_{\{l'\,|\,\boldsymbol{R}_{rl'}^{L}=1\wedge l'\preceq_{r}l\}}\boldsymbol{D}_{l'h}$.
$\boldsymbol{D}_{rn}^{aggr}$ codes the nominal (minimal) transport
time of route $r\in\mathbb{R}$ up to node $n\in\mathbb{N}$ and is
used to restrict the variable $x_{ntr}^{aggr}$ how much volume that
at most can have reached the node $n$ and thus the maximum volume
that may take the direct arcs, see constraint \textbf{Aggregate\_2}
below.
\end{description}

\subsubsection{\label{subsec:Variables}Variables}

The following variables are used in the model. The principal naming
schema is that $x$ is used for the flow layer , $y$ is used for
the routing layer  and $z$ is used for the demand layer. 
\begin{description}
\item [{$x_{rt}^{dep}$}] Matrix of real number values ranging over $\mathbb{R}\times\mathbb{T}$,
departure volumes at $t\in\mathbb{T}$ on route $r\in\mathbb{R}$
\item [{$x_{rt}^{arr}$}] Matrix of real number values ranging over $\mathbb{R}\times\mathbb{T}$,
arrival volumes at $t\in\mathbb{T}$ on route $r\in\mathbb{R}$
\item [{$x_{ntr}^{ext}$}] 3-dimensional array of real number values ranging
over $\mathbb{N}\times\mathbb{T}\times\mathbb{R}$, source and sink
for volumes taking route $r\in\mathbb{R}$ at time $t\in\mathbb{T}$
 $x^{ext}$works as a ``communication variable'' between demand
and implementation in terms of routes.
\item [{$x_{ltr}^{direct}$}] 3-dimensional array of real number values
ranging over $\mathbb{L}\times\mathbb{T}\times\mathbb{R}$.\\
These arcs corresponds to transport traversal on links within the
same time period holding the volume of transports. Note that the amount
can be fractional and that this is intentional, see section \ref{subsec:Traffic-volumes-and-capacity-usage}
for an explanation of the difference of volumes and individual trains.
\item [{$x_{ltr}^{next}$}] 3-dimensional array of real number values ranging
over $\mathbb{L}\times\mathbb{\widehat{T}}\times\mathbb{R}$.\\
These arcs corresponds to transport traversals on links from the current
time period to the next time period.
\item [{$x_{ntr}^{ni}$}] 3-dmensional array of real number values ranging
over $\mathbb{N}\times\mathbb{\widehat{T}}\times\mathbb{R}$.\\
These arcs corresponds to transports standing in a node from one time
period to the next. The superscript ``ni'' stands for ``node inventory''
\item [{$x_{ntr}^{in}$}] 3-dimensional array of real number  ranging over
$\mathbb{N}\times\mathbb{\widehat{T}}\times\mathbb{R}$.\\
These are additional redundant variables containing the sum of all
transportation arcs coming into a node. The superscript ``in'' stands
for ``into the node''.
\item [{$x_{ntr}^{aggr}$}] 3-dimensional array of real number values over
$\mathbb{N}\times\mathbb{\widehat{T}}\times\mathbb{R}$.\\
This is the aggregated maximum volumes of $r\in\mathbb{R}$ that could
possibly have reached the node $n\in\mathbb{N}$. $x^{aggr}$ is responsible
for pushing transports on to the $x^{next}$ arcs by restricting $x^{in}$
variables, otherwise it would be possible for the volumes to take
only direct arcs.
\item [{$y_{ot}^{post}$}] Matrix of real number values ranging over $\mathbb{D}\times\mathbb{\widehat{T}}$.\\
$y_{ot}^{post}$ holds the demanded volumes that are postponed to
depart later. Note that these arcs are similar to the $x^{ni}$ arc
on the same node for the relations, but $y_{ot}^{post}$ are arcs
on the route layer, not the implementation layer  (see section \ref{subsec:Demand-level},
\ref{subsec:Route-level} and \ref{subsec:Flow-level})
\item [{$y_{ot}^{cancel}$}] Matrix of canceled demands ranging over $\mathbb{D}\times\mathbb{T}$.\\
These variables hold the amount of cancellations of demands that cannot
be realized.
\item [{$y_{o}^{cancel}$}] Array of canceled demands ranging over $\mathbb{D}$.\\
$y_{o}^{cancel}$ holds the sum of all time frames in $y_{ot}^{cancel}$,
i.e. for a demand $o\in\mathbb{D}$ : $y_{o}^{cancel}=\sum_{t}y_{ot}^{cancel}$
\item [{$\boldsymbol{L}_{lth}^{C}$}] 3-dimensional matrix of real number
values ranging over $\mathbb{L}\times\mathbb{T}\times\mathbb{H}$.\\
This is the actual link capacity allocated to train type $h\in\mathbb{H}$
on link $l\in\mathbb{L}$ in period $t\in\mathbb{T}$, dependent on
the nominal capacity $\boldsymbol{C}_{lt}^{nom}$ and dependent on
dynamic restrictions on e.g. heterogeneous traffic, relation on direction
mix (for single track lines) and traffic restrictions due to track
work or other reasons. $\forall lt:\sum_{h}\boldsymbol{L}_{lth}^{C}\leq\boldsymbol{C}_{lt}^{nom}$.
\end{description}

\subsubsection{\label{subsec:Constraints}Constraints}

Constraints with names ``...alt...'' are alternatives, only one
should be used in the model.
\begin{description}
\item [{Capacity1}] Basic capacity constraint, the maximum nominal capacity
cannot be exceeded.
\[
\forall l\in\mathbb{L},t\in\mathbb{T}:(\sum_{h\in\mathbb{H}}\boldsymbol{L}_{lts}^{C})\leq\boldsymbol{C}_{lt}^{nom}
\]
\item [{Capacity2~alt1}] Alternative 1 of measuring capacity on single
track lines, where the volumes in opposite direction have to share
the capacity. 
\[
\forall l\in\mathbb{L},l<\sigma(l),t\in\mathbb{T}:\sum_{r\in\mathbb{R}\wedge\pi^{3}(\psi(r))=h}(\boldsymbol{L}_{lth}^{C}+L_{\sigma(l),th}^{C})\le0.5\,(\boldsymbol{C}_{lt}^{nom}+\boldsymbol{C}_{\sigma(l),t}^{nom})
\]
 
\item [{Capacity2~alt2}] Alternative 2 of measuring capacity on single
track lines, where the volumes in opposite direction have to share
the capacity. Add alternative constraint to \textbf{Capacity 1} above
for single line traffic:
\[
\begin{array}{l}
\forall l\in\mathbb{L},l<\sigma(l),t\in\mathbb{T}:\sum_{h\in H}\mathbf{L}_{lth}^{C}+\sum_{h\in H}\mathbf{L}_{\sigma(l)th}^{C}+w_{lt}\leq\boldsymbol{C}_{lt}^{nom}\\
\forall l\in\mathbb{L},l>\sigma(l),t\in\mathbb{T}:\sum_{h\in H}\mathbf{L}_{lth}^{C}+\sum_{h\in H}\mathbf{L}_{\sigma(l)th}^{C}+w_{\sigma(l)t}\leq\boldsymbol{C}_{lt}^{nom}\\
\forall l\in L,l>\sigma(l),t\in T:w_{lt}=0
\end{array}
\]
where $w_{lt}$ is the setup time for changing direction on the line.
$w_{lt}$ is given by the following constraints
\[
\begin{array}{l}
\forall l\in\mathbb{L},l<\sigma(l),t\in\mathbb{T}:0\leq\sum_{h\in\mathbb{H}}\mathbf{L}_{lth}^{C}\leq K_{lt}\,w_{lt}+M\,(1-\beta_{lt})\\
\forall l\in\mathbb{L},l<\sigma(l),t\in\mathbb{T}:0\leq\sum_{h\in\mathbb{H}}\mathbf{L}_{\sigma(l)th}^{C}\leq K_{\sigma(l)t}\,w_{lt}+M\,\beta_{lt}
\end{array}
\]
$\beta_{lt}$ is a new binary variable and $M$ a large constant
to dominate the equations so that exactly one of the two equations
is chosen. $K_{lt}\leq1$ is a setup coefficient, the smaller $K_{lt}$
the more does a change in direction affect the net capacity on the
link. It is important to note that $\beta_{lt}$ is a binary variable
making the problem a mixed integer problem. The only other variable
that could be argued should be integer is the demand cancel variables
$y^{cancel}$ depending on the interpretation of the variable if it
is interpreted as the probability of a cancellation (possibility to
accomplish the traffic) or if it is interpreted as an actual realization
of each required demand for traffic.
\item [{Capacity3}] Capacity on lines with heterogeneous traffic (i.e.
different speed), where resulting setup times between volumes with
different properties must be addressed.
\[
\forall l\in\mathbb{L},t\in\mathbb{T}:\sum_{h\in\mathbb{H}:p(l,h)}(\boldsymbol{L}_{lth}^{C}+\sum_{h'\in\mathbb{H}:h\ne h'\wedge p(l,h')}K\times\boldsymbol{L}_{lth'}^{C})\le\boldsymbol{C}_{lt}^{nom}
\]
where $p$ is a helper predicate, formulated as\\
$\forall l\in\mathbb{L},h\in\mathbb{H}:p(l,h)\equiv\forall r\in\mathbb{R}:(\boldsymbol{R}_{rl}^{L}=1\wedge\pi^{3}(\psi(r))=h)$\\
and $K$ is the ``heterogeneous coefficient'', tested with value
$K=0.25$. \\
The drawback of this formulation is that each $r\in\mathbb{R}$ will
have a setup time to each other $r'\in\mathbb{R}$, which is not the
case in reality where each change in performance due to different
train type $h$ and $h'$ (e.g. speed) gives rise to a setup time
if the corresponding train path are scheduled after each other. So,
for example, if there are 3 different train types, then there is a
sequence of them in the time period that will have two setup times
(and, possibly, two more to the adjacent time periods). But this formulation
will add capacity for 3{*}2=6 occurrences of performance change. With
increasing number of different train types the amount of added capacity
setup time will increase exponentially, which is not true in reality.
\item [{Capacity4}] Basic capacity restriction for each time frame and
track segment (arc) in the model. As an approximation, half of the
volume of the incoming ``next'' arcs counts as capacity consuming,
as well as half of the outgoing ``next'' arcs.
\[
\forall l\in\mathbb{L},t\in\mathbb{T},h\in\mathbb{H}:\sum_{r\in\mathbb{R}:\boldsymbol{R}_{rl}^{L}=1\land\pi^{3}(\psi(r))=h}(x_{ltr}^{direct}+0.5\,x_{l,t-1,r}^{next}+0.5\,x_{ltr}^{next})\leq\boldsymbol{L}_{lth}^{C}
\]
 
\item [{Demand1}] Initializing all postponed volumes at $t=0$ to $0.0$
i.e. if there are no postponed volumes entering the problem. 
\[
\forall d\in\mathbb{D}:y_{d,0}^{post}=0.0
\]
 
\item [{Demand2}] Initializing all demanded volumes in last time period
to zero.
\[
\forall d\in\mathbb{D}:y_{dt_{max}}^{post}=0.0
\]
By this, no demanded train can be postponed after the period end,
and thus it has to be canceled instead (giving rise to a large cost
by the objective function).
\end{description}
There is a possibility to initialize $y_{d,0}^{post}$ to non-zero
values if the problem formulation includes volumes that enters the
problem formulation from e.g. earlier time instances, It should also
be possible to identify $y_{dt_{max}}^{post}$ with $y_{d,0}^{post}$
which would then lead to a ``basic hourly pattern'' flow model (together
with other variables at the period start and end, e.g. $x_{n0r}^{ni}=x_{nt_{max}r}^{ni}$,
this has not been tested.
\begin{description}
\item [{Departure3}] Balance of each demand volume in each time frame at
demand origin nodes.
\[
\forall d\in\mathbb{D},t\in\mathbb{T}:\sum_{r\in\mathbb{I}(d)}(x_{rt}^{dep}-y_{d,t-1}^{post}+y_{dt}^{post}=\mathbf{D}_{dt}^{V}-y_{dt}^{cancel})
\]
\item [{Arrival1}] Requirement on total average runtime (duration) for
the volume for route $r\in\mathbb{R}$ in time period $t\in\mathbb{T}$
\[
\forall d\in\mathbb{D},t\in\mathbb{T},r\in\mathbb{R},n_{2}\in N:\pi^{2}(\psi(r))=n_{2}\wedge r\in\mathbb{I}(d)\rightarrow x_{rt}^{arr}\ge x_{n_{2},tr}^{aggr\_in}-x_{n_{2},t-1,r}^{aggr\_in}-S_{r}
\]
where $S_{r}$ is the allowed arrival slack. 
\item [{Cancel1}] The sum of all canceled volumes in each time frame is
the total volume of a certain demand $d\in\mathbb{D}$ being canceled
\[
\forall d\in\mathbb{D}:(\sum_{t\in\mathbb{T}}y_{dt}^{cancel})=y_{d}^{cancel}
\]
\item [{Cancel2}] Cancel route volumes, counted at the departure i.e. counting
canceled volumes at the source node
\[
\forall d\in\mathbb{D}:\sum_{t\in\mathbb{T}}(\sum_{r\in\mathbb{I}(d)\}}x_{rt}^{dep})=\nu(d)-y_{d}^{total}
\]
\item [{Cancel3}] Cancel route volumes, counted at the arrival i.e. counting
canceled volumes at the sink node
\[
\forall d\in\mathbb{D}:\sum_{t\in\mathbb{T}}(\sum_{r\in\mathbb{I}(d)\}}x_{rt}^{arr})=\nu(d)-y_{d}^{total}
\]
 This is a redundant constraint, as one of \textbf{Cancel2} and \textbf{Cancel3}
is sufficient.
\item [{Bound1}] All arcs that a route $r\in\mathbb{R}$ does not have
in its path are 0
\[
\forall l\in\mathbb{L},t\in\mathbb{T},r\in\mathbb{R},\boldsymbol{R}_{r,1}^{L}=0\}:x_{ltr}^{direct}=0
\]
\[
\forall l\in\mathbb{L},t\in\mathbb{T},r\in\mathbb{R},\boldsymbol{R}_{r,1}^{L}=0:x_{ltr}^{next}=0
\]
\item [{Bound2}] All node inventory arcs which are not possible for route
$r\in\mathbb{R}$ are 0
\[
\forall n\in\mathbb{N},t\in\mathbb{T},r\in\mathbb{R},\boldsymbol{R}_{rn}^{N}=0:x_{ntr}^{ni}=0
\]
 
\item [{Bound3}] For every route, all nodes except the source and sink
of the route can have volumes entering or exiting the flow layer.
\[
\forall n\in\mathbb{N},t\in\mathbb{T},d\in\mathbb{D},r\in\mathbb{R}:r\in\mathbb{I}(d)\wedge\pi^{1}(\psi(r))\ne n\wedge\pi^{2}(\psi(r))\ne n\rightarrow x_{ntr}^{ext}=0
\]
\item [{Bound4}] All next arcs starting at time zero (i.e. originates outside
the problem) are set to zero.
\[
\forall l\in\mathbb{L},r\in\mathbb{R}:x_{l,0,r}^{next}=0
\]
So are also all the next arcs passing out from the last time frame.
\[
\,\forall l\in\mathbb{L},r\in\mathbb{R}:x_{l,t_{max},r}^{next}=0
\]
 If these two, $\forall l\in\mathbb{L},r\in\mathbb{R}:x_{l,t_{max},r}^{next}=x_{l,0,r}^{next}$
are unified (together with some other flow variables) , then a ``basic
hourly pattern'' can be achieved.
\item [{Bound5}] All node inventory is zero at start.
\[
\forall n\in\mathbb{N},r\in\mathbb{R}:x_{n,0,r}^{ni}=0
\]
 If this is unified with $x_{l,t_{max},r}^{ni}$ (together with some
other flow variables) then a ``basic hourly pattern'' can be achieved
\item [{Flow1}] Source and sink correlating given demands and their possible
implementations in routes.
\[
\begin{array}{ll}
\forall n\in\mathbb{N},t\in\mathbb{T},r\in\mathbb{R}:x_{ntr}^{ext}=x_{rt}^{dep} & \text{ if }\pi^{1}(\psi(r))=n\\
\forall n\in\mathbb{N},t\in\mathbb{T},r\in\mathbb{R}:x_{ntr}^{ext}=-x_{rt}^{arr} & \text{ if }\pi^{2}(\psi(r))=n\\
\forall n\in\mathbb{N},t\in\mathbb{T},r\in\mathbb{R}:x_{ntr}^{ext}=0 & \text{ if }\pi^{1}(\psi(r))\ne n\wedge\pi^{2}(\psi(r))\ne n
\end{array}
\]
 Note that $\psi(r)$ is known as part of the problem statement and
therefore this case-construct can be dissolved.
\item [{Flow2}] Flow balance in each node
\[
\begin{array}{ll}
\forall n\in\mathbb{N},t\in\mathbb{T},r\in\mathbb{R},\boldsymbol{R}_{rn}^{N}=1:\,x_{ntr}^{ext} & \text{Flow in and out of flow layer (source and sink)}\\
+x_{n,t-1,r}^{ni} & \text{Node inventory from previous period}\\
+\sum_{l\in\mathbb{L}\wedge h(l)=n}(x_{ltr}^{direct}+x_{l,t-1,r}^{next}) & \text{Flow on oncoming line arcs into this node}\\
-x_{ntr}^{ni} & \text{Node inventory gouing to next period}\\
-\sum_{l\in\mathbb{L}\wedge t(l)=n\}}(x_{ltr}^{direct}+x_{ltr}^{next})=0 & \text{Flow on outgoing line arcs from this node}
\end{array}
\]
\item [{Flow3}] Share of volumes arriving at each node in each time period
\[
\forall n\in\mathbb{N},t\in\mathbb{T},r\in\mathbb{R},\pi^{1}(\psi(r))=n:x_{ntr}^{in}=x_{rt}^{dep}+\sum_{l\in\mathbb{L}\wedge t(l)=n\wedge\boldsymbol{R}_{rl}^{L}=1}(x_{i,t-1,r}^{next}+x_{itr}^{direct})
\]
\[
\forall n\in\mathbb{N},t\in\mathbb{T},r\in\mathbb{R},\pi^{1}(\psi(r))\ne n:x_{ntr}^{in}=\sum_{l\in\mathbb{L}\wedge t(l)=n\wedge\boldsymbol{R}_{ri}^{L}=1}(x_{i,t-1,r}^{next}+x_{itr}^{direct})
\]
\item [{Aggregate1}] At time 0 no trains have left their origin and thus
no train can have reached any other station
\[
\forall n\in\mathbb{N},r\in\mathbb{R}:x_{n,0,r}^{aggr}=0
\]
 
\item [{Aggregate2.1}] \begin{flushleft}
All nodes $n$ where route $r$ does not pass are
set to 0
\[
\forall n\in\mathbb{N},t\in\mathbb{T},r\in\mathbb{R},\boldsymbol{R}_{rn}^{N}=0:x_{ntr}^{aggr}=0
\]
\par\end{flushleft}
\item [{Aggregate2.2}] \begin{flushleft}
Aggregated sums along the route's path. For each
$r\in\mathbb{R},$ node $n\in\mathbb{N}$ in $r$'s path and time
period $t\in\mathbb{T}$ calculate how much volume that at most can
have reached the node along $r$
\[
\begin{array}{ll}
\forall n\in\mathbb{N},t\in\mathbb{T},r\in\mathbb{R},\boldsymbol{R}_{rn}^{N}=1:x_{ntr}^{aggr}= & \text{aggregation at time\ensuremath{t}}\\
\,x_{n,t-1,r}^{aggr} & \text{all previous aggregated}\\
 & \text{volumes}\\
+\,x_{rt}^{dep}\times max(1-\boldsymbol{D}_{rn}^{aggr},0) & \text{Volumes departed in the}\\
 & \text{same time period}\\
+\,\sum_{0<t'<t\wedge t-t'=\underline{\boldsymbol{D}_{rn}^{aggr}}}x_{r,t'}^{dep}\times(\overline{\boldsymbol{D}_{rt}^{aggr}}-\boldsymbol{D}_{rn}^{aggr}) & \text{Volumes departed in earlier}\\
 & \text{periods, earlier part}\\
+\,\sum_{0<t'<t\wedge t-t'=\overline{\boldsymbol{D}_{rn}^{aggr}}}x_{r,t'}^{dep}\times(\boldsymbol{D}_{rn}^{aggr}-\underline{\boldsymbol{D}_{rn}^{aggr}}) & \text{Volumes departed in earlier}\\
 & \text{periods, later part}
\end{array}
\]
where $\overline{\boldsymbol{D}_{rn}^{aggr}}$ is the integer ceiling
of $\boldsymbol{D}_{rn}^{aggr}$ and $\underline{\boldsymbol{D}_{rn}^{aggr}}$
is the integer floor of $\boldsymbol{D}_{rn}^{aggr}$. Note that $\boldsymbol{D}_{rn}^{aggr}$
is already known when the problem is initiated, as are the $max(...)$
term.
\par\end{flushleft}
\item [{Aggregate3}] Restriction on volumes that may take the direct arc
$x_{direct}$, first time period.
\[
\forall n\in\mathbb{N},r\in\mathbb{R},R_{rn}^{N}=1:x_{n,1,r}^{aggr}\geq x_{n,1,r}^{in}
\]
 The constraint works in conjunction with the \textbf{Flow3} constraint
to pose the restriction on the $x^{direct}$ arcs.
\item [{Aggregate4}] Restriction on volumes that may take the direct arcs
$x^{direct}$ , all but first time period.
\[
\forall n\in\mathbb{N},t\in\mathbb{T},r\in\mathbb{R},R_{rn}^{N}=1:t>1\rightarrow x_{ntr}^{aggr}\geq\sum_{t'\in\mathbb{T}:t'\leq t}x_{n,t',r}^{in}
\]
 The constraint works in conjunction with the \textbf{Flow3} constraint
to pose the restriction on the $x^{direct}$ arcs.
\end{description}

\subsubsection{\label{subsec:Objective-function}Objective function }

The objective function should reflect the efficiency of the traffic
and preferably also the benefit for the owner of the infrastructure
and trains. To cancel trains should be the ultimate offer if the capacity
is not enough, and therefore the cancellation coefficient is chosen
high. Postponing departure compared to the the applied departure time
also comes with a cost, as well as the running times regardless of
the route taken. It is also possible to add terms for e.g. distance
traveled (known at the routing level and a constant for each described).

The following objective functions has been used during the development
of the model and is used in the first example in section \ref{sec:Examples}.

\[
obj=\sum_{d\in\mathbb{D}}(1000\times C_{d}^{total}+\sum_{t\in\mathbb{\widehat{T}}}20\times y_{dt}^{post})+\sum_{d\in\mathbb{D}}(\sum_{r\in\mathbb{I}(d)}(\sum_{t\in\mathbb{T}}t\times x_{rt}^{arr}-\sum_{t\in\mathbb{T}}t\times x_{rt}^{dep}))/\sum_{d\in\mathbb{D}}\nu(d)
\]

For the second example showing one of the intended use of the model
a variant of the priority categories, used by Trafikverket \cite{Trafikverket_JNB2021},
to solve disputes when capacity conflicts cannot be resolved by negotiations,
is used. These categories are based on a socio-economic framework
call ASEK \cite{Trafikverket_ASKE7.0} and give standard values for
departure delays, extended running times and train cancellations,
so the framework, or a variant of it, can be used to form an objective
function. 

\section{Complexity}

The complete model proposed in this paper can at a first glance look
complex, not only to understand but also from a complexity point of
view. However, the crucial step to control the complexity is the restriction
at the flow level network to only use named routes given by the Train
Service Classes, TCSs. With these routes the choice of path through
the geographical network is given, and what is left on the flow level,
per route, is to schedule the flow onto the three different arc types
direct arcs, next arcs and node inventory arcs, per link. The choice
of route is made among the given named routes, where in practice only
a few are relevant to state.

Also the two aggregation sums restrict the solution space. The maximum
aggregated sum at each node is crucial to have volumes obey the timing
constraints of their train types with tightening bounds. The same
is true for the minimum aggregated sums as well, i.e. the constraint
at each node that states that at least a certain volume of trains
has reached the node at a certain point in time.

The number of direct and next arcs are dependent on the number of
time periods, the number of network links in the railway network and
the number of routes using the network link. The number of node inventory
variables as well as the aggregated sum variables are dependent on
the number of network nodes, time periods and routes passing the node.
For all these four variable types the number of routes that passes
the link or node is restricted to the routes relevant to use for a
demand, and can therefore be assumed to be restricted to the routes
that passes each node and link. This means that the number of time
periods and network links/nodes are the dominant factors for the number
of variables. All these variables are also real number variables,
i.e. there is no integral property requirement on them.

At the demand level the cancellations variables (of whole trains)
must be integer declared. Note that it is possible that the postpone
variables does not have to be integers, since a fractional postponing
value means that the train is postponed in time a little bit and not
a complete time period. There is however an uninvestigated relation
between volumes postponed over several time periods, if this should
be allowed or not. It would perhaps be natural to restrict postponing
parts of volumes of trains only one time period, at least not let
the volume be spread over several time periods with holes in departure
volumes in the middle. One initial idea is that postponed volumes
to the next time period is dependent on the added volumes in the same
period together with the volumes departing in this time period. This
is currently not solved and is part of future work.

There are a number of constraint groups at the flow level that in
principle could grow large, see section \ref{sec:Mathematical-model}.
As for variables, the main source for growing numbers are the number
of nodes and links together with the number of time periods. Again,
by restricting the number of routes to those that are relevant to
investigate the number of constraints are restricted to the ones that
are relevant for the problem at hand. 

\section{Examples\label{sec:Examples}}

We present two examples, one smaller untended to illustrate the model,
and one larger to illustrate one of the intended uses of it. 

\subsection{Illustration of the model\label{subsec:Illustration-of-the}}

\begin{figure}[H]
\begin{minipage}[c]{0.4\columnwidth}%
\includegraphics[scale=0.25]{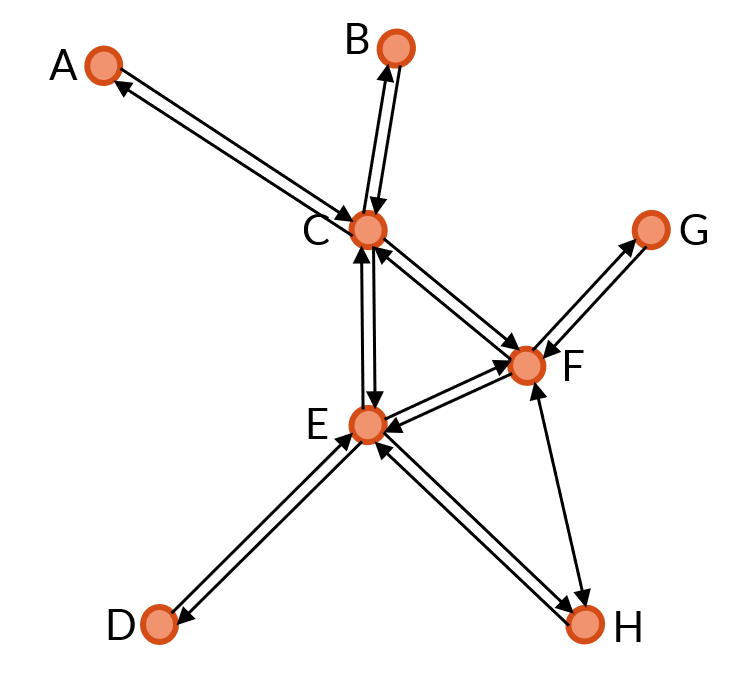}%
\end{minipage}~~~%
\begin{tabular}{|c|c|c|c|c|c|c|c|}
\hline 
 & 1 & 2 & 3 & 4 & 5 & 6 & 7\tabularnewline
\hline 
\hline 
A-C & 5 & 5 & 5 & 5 & 5 & 5 & 5\tabularnewline
\hline 
C-A & 5 & 5 & 5 & 5 & 5 & 5 & 5\tabularnewline
\hline 
B-C & 5 & 5 & 5 & 5 & 5 & 5 & 5\tabularnewline
\hline 
C-B & 5 & 5 & 5 & 5 & 5 & 5 & 5\tabularnewline
\hline 
C-E & 5 & 5 & 5 & 5 & 5 & 5 & 5\tabularnewline
\hline 
E-C & 5 & 5 & 5 & 5 & 5 & 5 & 5\tabularnewline
\hline 
C-F & 5 & 5 & 5 & 5 & 5 & 5 & 5\tabularnewline
\hline 
F-C & 5 & 5 & 5 & 5 & 5 & 5 & 5\tabularnewline
\hline 
D-E & 5 & 5 & 5 & 5 & 5 & 5 & 5\tabularnewline
\hline 
E-D & 5 & 5 & 5 & 5 & 5 & 5 & 5\tabularnewline
\hline 
F-G & 5 & 5 & 5 & 5 & 5 & 5 & 5\tabularnewline
\hline 
G-F & 5 & 5 & 5 & 5 & 5 & 5 & 5\tabularnewline
\hline 
E-F & 5 & 5 & 5 & 5 & 5 & 5 & 5\tabularnewline
\hline 
F-E & 5 & 5 & 5 & 5 & 5 & 5 & 5\tabularnewline
\hline 
E-H & 5 & 5 & 5 & 5 & 5 & 5 & 5\tabularnewline
\hline 
H-E & 5 & 5 & 5 & 5 & 5 & 5 & 5\tabularnewline
\hline 
F-H & \multirow{2}{*}{5} & \multirow{2}{*}{5} & \multirow{2}{*}{5} & \multirow{2}{*}{5} & \multirow{2}{*}{5} & \multirow{2}{*}{5} & \multirow{2}{*}{5}\tabularnewline
\cline{1-1} 
H-F &  &  &  &  &  &  & \tabularnewline
\hline 
\end{tabular}

\caption{\label{fig:Example-basic-railway}Example basic railway network}
\end{figure}

The following small example illustrates the model. It is chosen to
be small, yet illustrate the different properties of a railway network
and opportunities to vary traffic solutions. In figure \ref{fig:Example-basic-railway}
the basic network is shown, with the basic characteristics given in
the table. Note that there exists one single track line, arcs F-H
and H-F use the same track, which is shown in the table where the
two directions share the same capacity. This capacity is then dependent
on the mix of directions, some will be used for setup times to change
direction.

\begin{figure*}[htpb]

\begin{minipage}[t][1\totalheight][b]{0.45\columnwidth}%
\begin{figure}[H]
\includegraphics[scale=0.25]{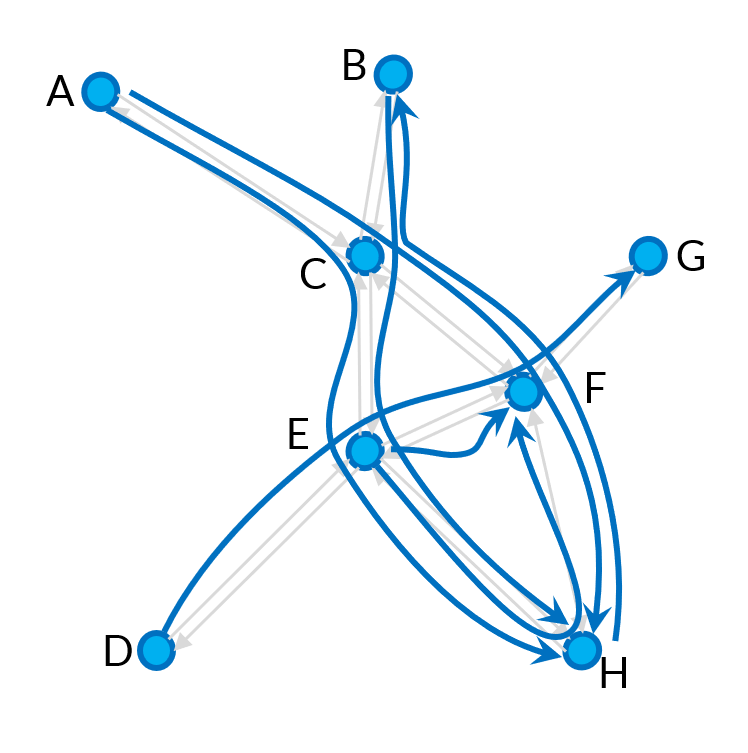}

\caption{\label{fig:Example-basic-railway-1}Example basic railway network}
\label{Example routes}
\end{figure}
\end{minipage}~~~%
\begin{minipage}[t]{0.45\columnwidth}%
\begin{figure}[H]
\includegraphics[scale=0.25]{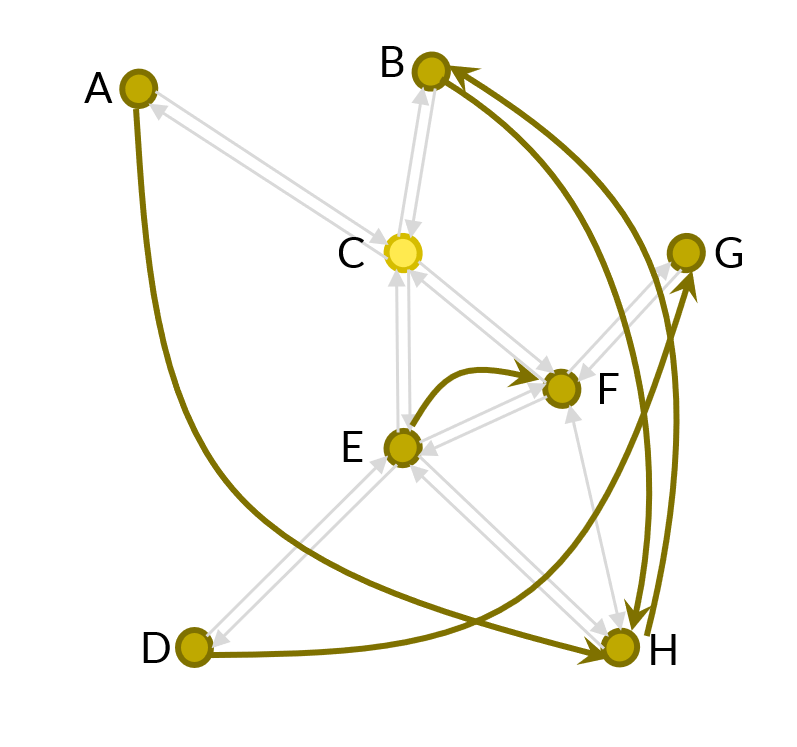}

\caption{\label{fig:Example-demands}Example demands}
\label{Example demand 1}
\end{figure}
\end{minipage}

\end{figure*}

The demands and routes are given in figure \ref{Example routes} and
figure \ref{Example demand 1} with data given in tables \ref{Table Example Routes}
and \ref{Table Example Demands 1} respectively. The node 'C' in figure
\ref{Table Example Demands 1} that is lighter in colour is not part
of the demand level, since no traffic originates or ends at this node.
In a real example much more of the traffic should be of the kind that
goes back and forth between the nodes B and H, reflecting a regional
traffic scenario. The purpose of this example is however to illustrate
the model rather than having all setup as close to reality as possible
everywhere.

\begin{table}[H]
\begin{tabular}{|>{\centering}p{1.3cm}|>{\centering}p{0.2cm}|>{\centering}p{0.2cm}|>{\centering}p{0.2cm}|>{\centering}p{0.2cm}|>{\centering}p{0.2cm}|>{\centering}p{0.2cm}|>{\centering}p{0.2cm}|>{\centering}p{0.2cm}|>{\centering}p{0.2cm}|>{\centering}p{0.2cm}|>{\centering}p{0.2cm}|}
\hline 
 & \begin{turn}{90}
A-C
\end{turn} & \begin{turn}{90}
B-C
\end{turn} & \begin{turn}{90}
C-E
\end{turn} & \begin{turn}{90}
C-F
\end{turn} & \begin{turn}{90}
E-H
\end{turn} & \begin{turn}{90}
F-H
\end{turn} & \begin{turn}{90}
D-E
\end{turn} & \begin{turn}{90}
E-F
\end{turn} & \begin{turn}{90}
F-G
\end{turn} & \begin{turn}{90}
H-F
\end{turn} & \begin{turn}{90}
F-C
\end{turn}\tabularnewline
\hline 
\hline 
A-H-f1 & 1 & 0 & 1 & 0 & 1 & 0 & 0 & 0 & 0 & 0 & 0\tabularnewline
\hline 
A-H-f2 & 1 & 0 & 0 & 1 & 0 & 1 & 0 & 0 & 0 & 0 & 0\tabularnewline
\hline 
B-H-p1 & 0 & 1 & 1 & 0 & 1 & 0 & 0 & 0 & 0 & 0 & 0\tabularnewline
\hline 
H-B-p1 & 0 & 0 & 0 & 0 & 0 & 0 & 0 & 0 & 0 & 1 & 1\tabularnewline
\hline 
D-G-f1 & 0 & 0 & 0 & 0 & 0 & 0 & 1 & 1 & 1 & 0 & 0\tabularnewline
\hline 
E-F-p1 & 0 & 0 & 0 & 0 & 0 & 0 & 0 & 1 & 0 & 0 & 0\tabularnewline
\hline 
E-F-p2 & 0 & 0 & 0 & 0 & 1 & 0 & 0 & 0 & 0 & 1 & 0\tabularnewline
\hline 
\end{tabular}\caption{\label{tab:Example,-possible-(named)}Example, possible (named) routes
to implement the demands}
\label{Table Example Routes}
\end{table}

\begin{table}[H]
\begin{tabular}{|>{\centering}p{1.2cm}|>{\centering}p{0.26cm}|>{\centering}p{0.26cm}|>{\centering}p{0.26cm}|>{\centering}p{0.26cm}|>{\centering}p{0.26cm}|>{\centering}p{0.26cm}|>{\centering}p{0.26cm}|}
\hline 
 & 1 & 2 & 3 & 4 & 5 & 6 & 7\tabularnewline
\hline 
\hline 
A-H-f & 2 & 1 & 0 & 0 & 1 & 2 & 0\tabularnewline
\hline 
B-H-p & 0 & 1 & 3 & 2 & 3 & 1 & 0\tabularnewline
\hline 
H-B-p & 0 & 1 & 3 & 2 & 3 & 1 & 0\tabularnewline
\hline 
D-G-f & 3 & 0 & 0 & 0 & 1 & 3 & 0\tabularnewline
\hline 
E-F-p & 0 & 0 & 1 & 1 & 1 & 0 & 0\tabularnewline
\hline 
\end{tabular}

\caption{\label{tab:Example,-demands-for}Example, demands for traffic}
\label{Table Example Demands 1}
\end{table}

Running this small example in the modeling framework Minizinc \cite{Minizinc}
with the solver CBC \cite{Coin-OR} takes less than 1 second, and
the result (capacity usage) is shown in table \ref{tab:Result-from-executing}.
Empty rows does not get any traffic. The objective value is 0.4694.

\begin{table}
\begin{tabular}{|c|c|c|c|c|c|c|c|}
\hline 
 & 1 & 2 & 3 & 4 & 5 & 6 & 7\tabularnewline
\hline 
\hline 
A-C & 1,7 & 1,15 & 0,15 & 0 & 0,85 & 1,2 & 0,95\tabularnewline
\hline 
C-A &  &  &  &  &  &  & \tabularnewline
\hline 
B-C & 0 & 0,95 & 2,15 & 2,8 & 2,95 & 1 & 0,15\tabularnewline
\hline 
C-B & 0 & 0,7 & 2,4 & 2,4 & 2,6 & 1,6 & 0,3\tabularnewline
\hline 
C-E & 0 & 0,8 & 2,3 & 2,6 & 2,7 & 1,3 & 0,3\tabularnewline
\hline 
E-C &  &  &  &  &  &  & \tabularnewline
\hline 
C-F & 1,05 & 1,48 & 0,48 & 0 & 0,52 & 0,58 & 1,9\tabularnewline
\hline 
F-C & 0 & 0,8 & 2,3 & 2,5 & 2,8 & 1,3 & 0,3\tabularnewline
\hline 
D-E & 2,55 & 0,45 & 0 & 0 & 0,65 & 2,3 & 1.05\tabularnewline
\hline 
E-D &  &  &  &  &  &  & \tabularnewline
\hline 
F-G & 1,2 & 1,8 & 0 & 0 & 0,3 & 1,6 & 2,1\tabularnewline
\hline 
G-F &  &  &  &  &  &  & \tabularnewline
\hline 
E-F & 1,8 & 1,2 & 0,95 & 1 & 1,25 & 1,7 & 2,1\tabularnewline
\hline 
F-E &  &  &  &  &  &  & \tabularnewline
\hline 
E-H & 0 & 0,7 & 2,4 & 2,4 & 2,6 & 1,6 & 0,3\tabularnewline
\hline 
H-E &  &  &  &  &  &  & \tabularnewline
\hline 
F-H & 0,4 & 1,8 & 0,8 & 0 & 0,2 & 0,9 & 1,9\tabularnewline
\hline 
H-F & 0 & 0,95 & 2,6 & 2,05 & 3,25 & 0,85 & 0,3\tabularnewline
\hline 
\end{tabular}

\caption{\label{tab:Result-from-executing}Result from executing the model
on the example}
\end{table}

If we do a slight change and set the capacity availability on link
E-F to 0 during time period 4, simulating e-g- a TCR during that period,
we get a new solution with slighter higher objective of 0.48. Since
there cannot be any volumes taking the route E-F any longer, the demanded
volumes of the traffic E-F previously taking the route E-F-p1 now
partially is taking the route E-F-p2 instead which already had volumes
for the route E-H-F. If there would have been volumes demanded for
D-G-f they would have had to wait somewhere before E (or run slowly)
since there is only one named route D-G-f1 for D-G-f. In other words,
no named rerouting possibilities have been given since no other route
is present in the input data of the example. 

\begin{table}[H]
\begin{minipage}[c][1\totalheight][t]{0.45\columnwidth}%
\includegraphics[scale=0.3]{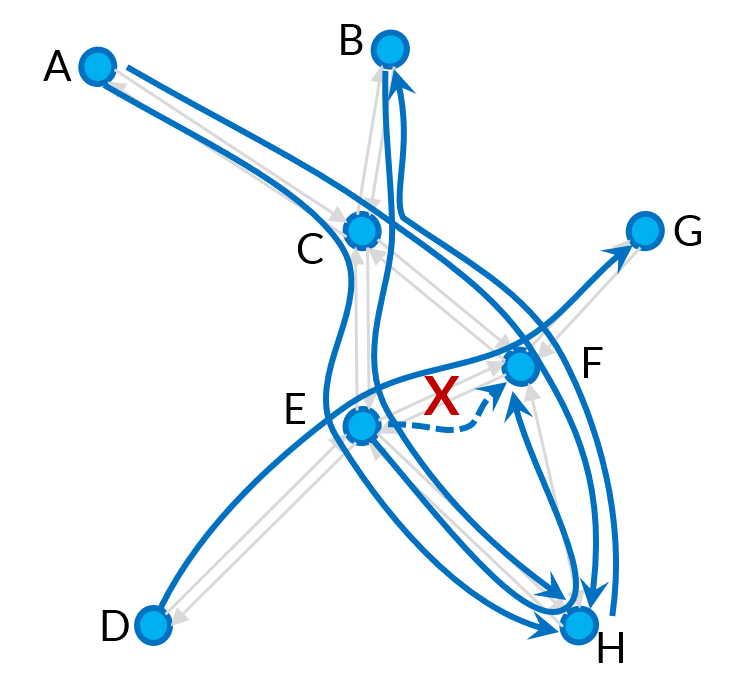}%
\end{minipage}%
\begin{tabular}{|>{\centering}p{0.7cm}|>{\centering}p{0.55cm}|>{\centering}p{0.55cm}|>{\centering}p{0.55cm}|>{\centering}p{0.55cm}|>{\centering}p{0.55cm}|>{\centering}p{0.55cm}|>{\centering}p{0.55cm}|}
\hline 
 & 1 & 2 & 3 & 4 & 5 & 6 & 7\tabularnewline
\hline 
\hline 
A-C & 1,7 & 1,15 & 0,15 & 0 & 0,85 & 1,2 & 0,95\tabularnewline
\hline 
C-A &  &  &  &  &  &  & \tabularnewline
\hline 
B-C & 0 & 0,9 & 2,65 & 2,35 & 2,95 & 1 & 0,15\tabularnewline
\hline 
C-B & 0 & 0,7 & 2,4 & 2,3 & 2,7 & 1,65 & 0,25\tabularnewline
\hline 
C-E & 0 & 0,8 & 2,3 & 2,6 & 2,7 & 1,3 & 0,3\tabularnewline
\hline 
E-C &  &  &  &  &  &  & \tabularnewline
\hline 
C-F & 1,05 & 1,33 & 0,63 & 0 & 0,38 & 0,72 & 1,9\tabularnewline
\hline 
F-C & 0 & 0,8 & 2,6 & 2 & 2,7 & 1,7 & 0,2\tabularnewline
\hline 
D-E & 2,55 & 0,45 & 0 & 0 & 0,85 & 2,1 & 1,05\tabularnewline
\hline 
E-D &  &  &  &  &  &  & \tabularnewline
\hline 
F-G & 0,9 & 2,1 & 0 & 0 & 0,4 & 1,5 & 2,1\tabularnewline
\hline 
G-F &  &  &  &  &  &  & \tabularnewline
\hline 
E-F & 1,5 & 1,5 & 0,9 & \textcolor{red}{0} & 1,65 & 1,35 & 2,1\tabularnewline
\hline 
F-E &  &  &  &  &  &  & \tabularnewline
\hline 
E-H & 0 & 0,7 & 2,4 & 3,2 & 2,8 & 1,6 & 0,3\tabularnewline
\hline 
H-E &  &  &  &  &  &  & \tabularnewline
\hline 
F-H & 0,4 & 1,65 & 0,95 & 0 & 0,05 & 1,05 & 1,9\tabularnewline
\hline 
H-F & 0 & 0,9 & 2,8 & 2,9 & 2,95 & 1,35 & 0,1\tabularnewline
\hline 
setup & 0 & 0,9 & 0,95 & 0 & 0,05 & 1,05 & 0,1~\tabularnewline
\hline 
\end{tabular}

\caption{\label{tab:Results-from-executing}Results from executing the model
with capacity availability in time period 4 to 0 on link E-F}
\end{table}

Note that the setup time for changing direction is included in the
last row of the table shown in table \ref{tab:Results-from-executing}.
The example uses the \textbf{Capacity2\_alt2} constraint, with the
setup coefficient $C_{lt}=1$, therefore the setup time will equal
the lower volume value of the two directions. The total capacity used
on a single track line is the sum of the volume in both directions
plus the setup time, and this sum cannot be larger than the capacity
stated, for each time period in the problem.

\subsection{Example use of the model\label{subsec:Example-use-of}}

\begin{figure}
\includegraphics[scale=0.35]{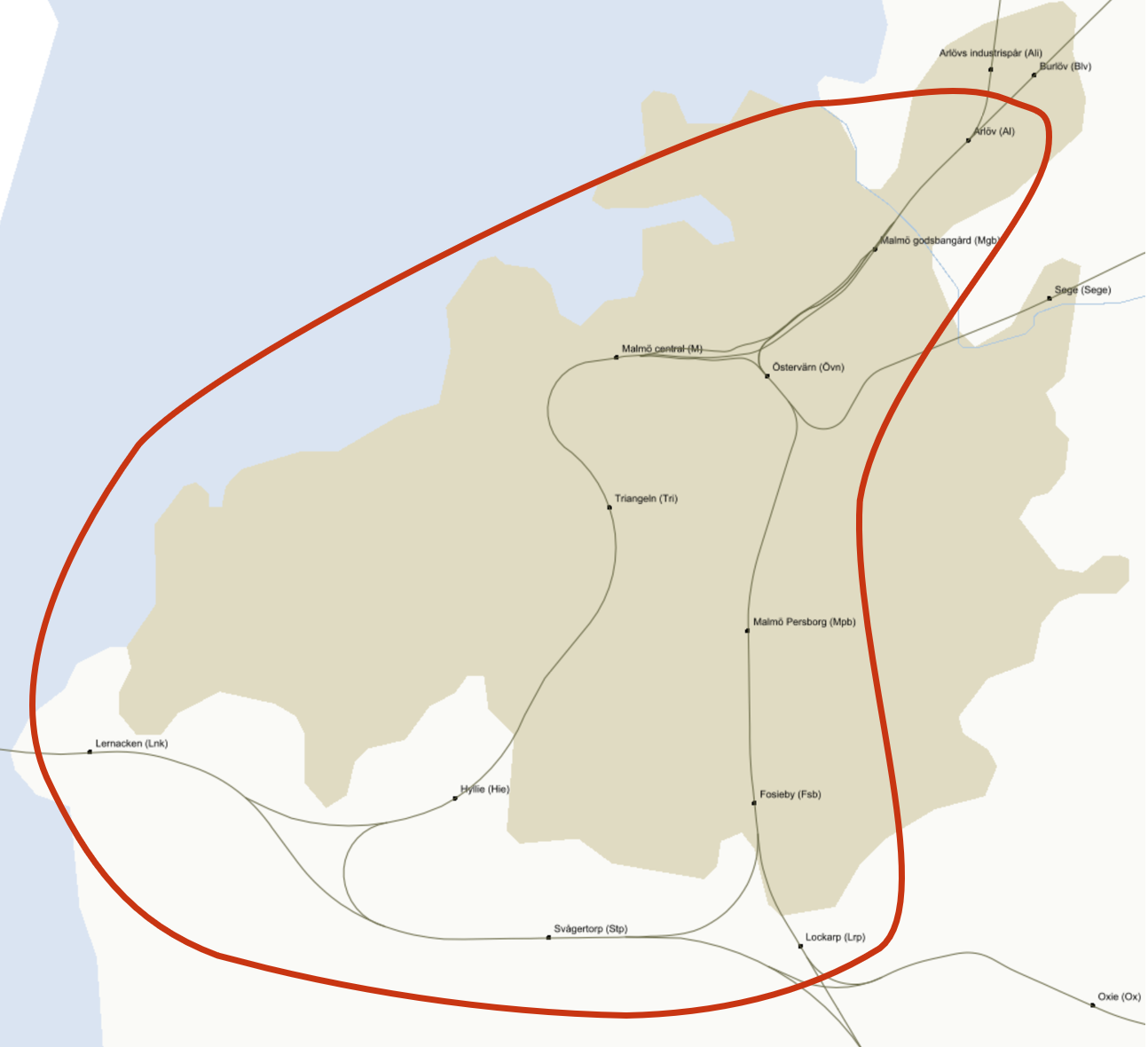}

\caption{\label{fig:Region-around-Malm=0000F6}Region around Malm�}
\end{figure}

The following example illustrates the use of the model in a real network.
The example uses a region around Malm� in Sweden, see figure \ref{fig:Region-around-Malm=0000F6}
where the basic network is shown. 

The example uses two days of traffic, discretizised into 120 minute
periods, and the aim is to see the effect on the traffic over the
bridge Malm� (Lernacken, LNK) to Copenhagen (Pepparholmen, PHM, last
timetable point in Sweden) when restricting the double track line
HIE - M to single track line. The first 6 periods are only to initialize
the ``pressure'' on the system and analogously the last 6 periods
are only for ``chilling down'' all the flows on the different arcs
in the network. The middle 12 time periods consists of the actual
measure of the impact when restricting the capacity.

The problem is based on the timetable 2020's traffic and consists
of 13 line segments giving 26 directed arcs in the network. There
are three different classes of trains: freight trains (FT), intercity
trains (IC) and regional trains (RT). In total there are 28 different
service types demanded, based on these classes and the origin-destination
pairs including different stop patterns of the trains. In total there
are 1136 individual demanded trains. In the basic scenario the routes
are determined, which gives us 28 TCSs implementing all the demands.

Solving this with the solver CBC, an open source MIP solver, takes
about 10 seconds on a standard laptop. No cancellations happen which
is expected since this was actual traffic originally planned in the
timetable without restricting the capacity.

\begin{figure}[H]
\includegraphics[scale=0.25]{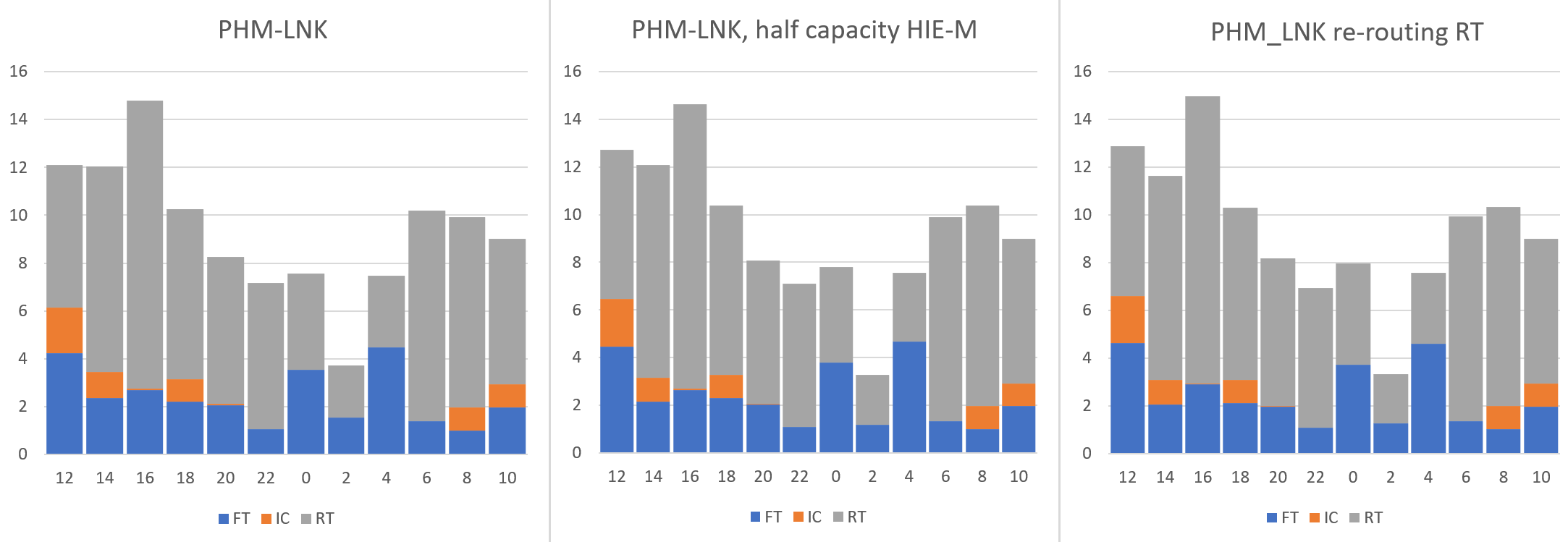}

\caption{Capacity usage on the link PHM-LNK for the three different train classes
FT, IC and RT\label{fig:Capacity-usage-on-PHM-LNK}}
\end{figure}

\begin{figure}[H]
\includegraphics[scale=0.25]{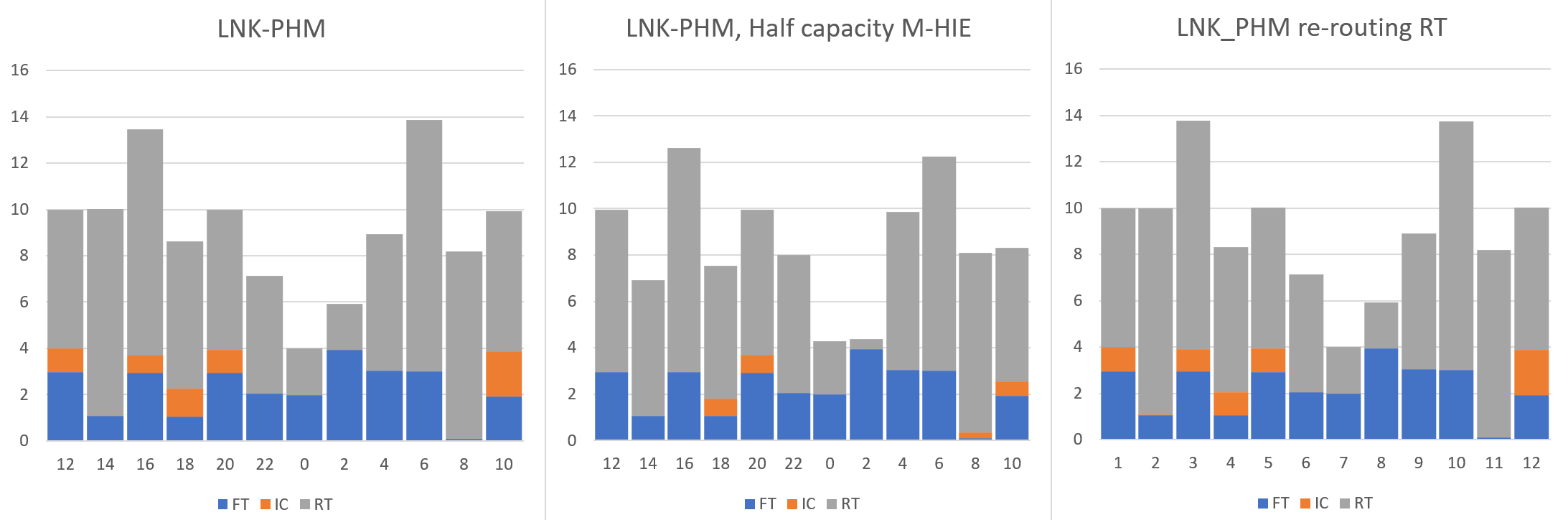}

\caption{Capacity usage on the link LNK-PHM for the three different train classes
FT, IC and RT\label{fig:Capacity-usage-on-LNK-PHM}}
\end{figure}

Now we restrict the capacity to half in both directions on the link
HIE-M. Without any rerouting of any trains, 19 trains have to be canceled
in total. However, if a rerouting possibility is introduced for the
regional trains going PHM-AL and AL-PHM to go round Malm� instead
(which is, however, not a preferred route since it is both longer
and it involves a direction change at Malm� central station) then
all trains in the original traffic can be accomplished but sometimes
with longer journey times. Figures \ref{fig:Capacity-usage-on-PHM-LNK}
and \ref{fig:Capacity-usage-on-LNK-PHM} shows the capacity usage
over time for the link PHM-LNK and LNK-PHM for the different example
runs.

The changes are mostly seen on the link LNK-PHM just after midnight.

This is one intended use of the model, to investigate and do scenario
analysis of different rerouting possibilities and measure the impact
of them. Experiments have been performed with both scenarios with
possessions and restrictions in capacity as well as introductory tests
within the Minimum Viable Product 6 in TTR \cite{TTR} for the case
Denmark-Sweden-Norway with the aim to facilitate the construction
of a capacity model and capacity supply.

\section{\label{sec:Time-period-length}Discussion on the time period length}

\subsection{Extending travel over more than one period}

The time period length should is crucial to choose. The smaller it
is, the better resolution the model will have. Larger time periods
decreases the complexity of the model. A requirement on the time period
as the model is formulated in section \ref{sec:Mathematical-model}
is that no volume's traversal time in the investigated network can
be larger than the time period length (i.e. it must reach the next
node within one time frame):

$\forall lt:x_{lt}^{ut}\leq c_{lt}^{ut}$

\begin{figure}[H]
\begin{centering}
\includegraphics[scale=0.3]{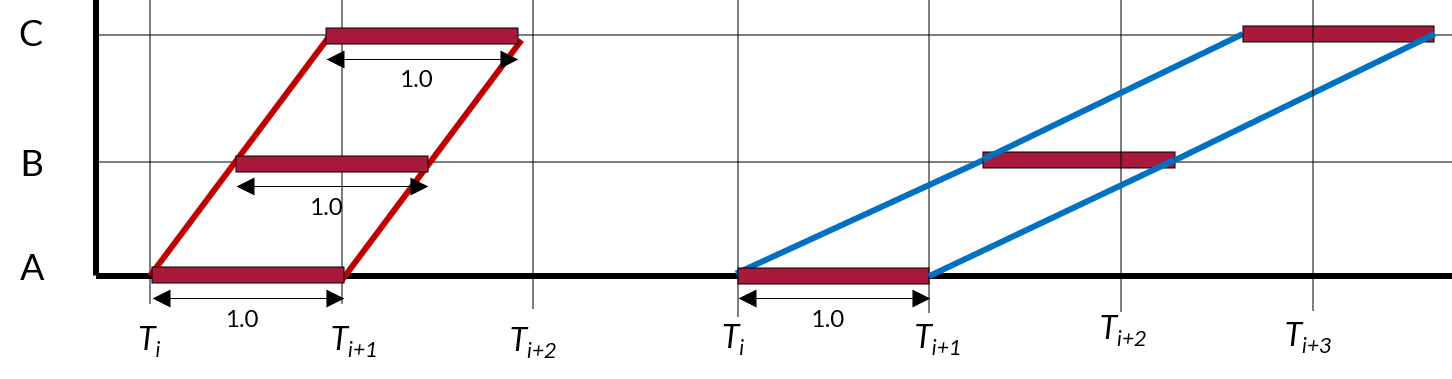}
\par\end{centering}
\caption{\label{fig:Time-period-length}Time period length and train's traversal
time}
\end{figure}

The reason is that there are no arcs ``jumping'' or ``skipping''
the next time period. As can be seen in the figures of section \ref{subsec:Traffic-volumes-and-capacity-usage}
the flow is distributed over the direct and next arc, but there are
no arc skipping one time period, as would be necessary to have a time
period smaller than the traversal time of the modeled trains (again
we omit the node inventory arcs, since they only add the possibility
to stand still at the station). The difference is highlighted in figure
\ref{fig:Time-period-length}.

The model can be extended with additional arcs skipping one or more
time periods, shown in figure \ref{fig:Time-period-length,}. If this
is done, the constraints in the model have to be updated with this
extension and the model gets more complex. For instance the capacity
constraints has to be extended to cope with all arcs that ``uses''
(or ``touches'') the time period, including the arcs that ``flies''
over the time period (the trains do travel in the time period, although
the ``jumping'' arc does not start or end within the time period).

\begin{figure}[H]
\begin{centering}
\includegraphics[scale=0.35]{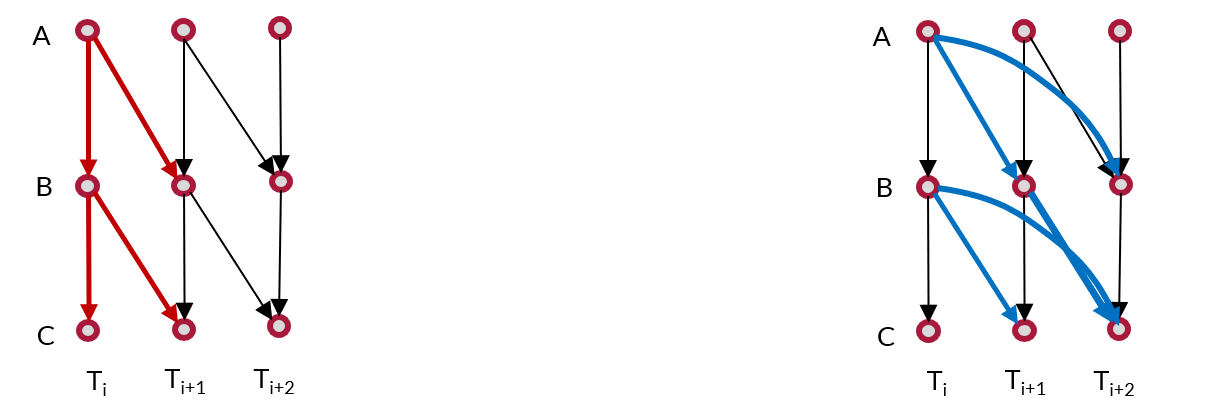}
\par\end{centering}
\caption{\label{fig:Time-period-length,}Time period length, ``jumping arcs''
for slow train paths}
\end{figure}

\subsection{Large departure time windows}

\begin{figure}[H]
\begin{centering}
\includegraphics[scale=0.3]{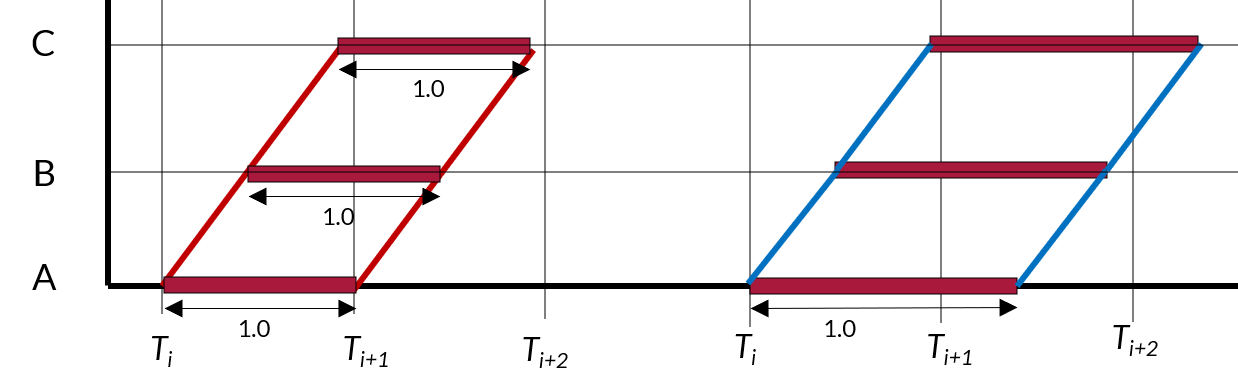}
\par\end{centering}
\caption{\label{fig:Time-windows-larger}Time windows larger than time period}
\end{figure}

Volumes for trains with large departure time windows, larger than
one time period, can be modeled to stretch over more than one time
period. It is just to distribute the volume over more than one time
period. Note however that when doing so in the current model the aggregated
volumes must be calculated with the knowledge that the stretch of
the volume has a certain form, see figure \ref{fig:Time-windows-larger}.

\begin{figure}[H]
\begin{centering}
\includegraphics[scale=0.3]{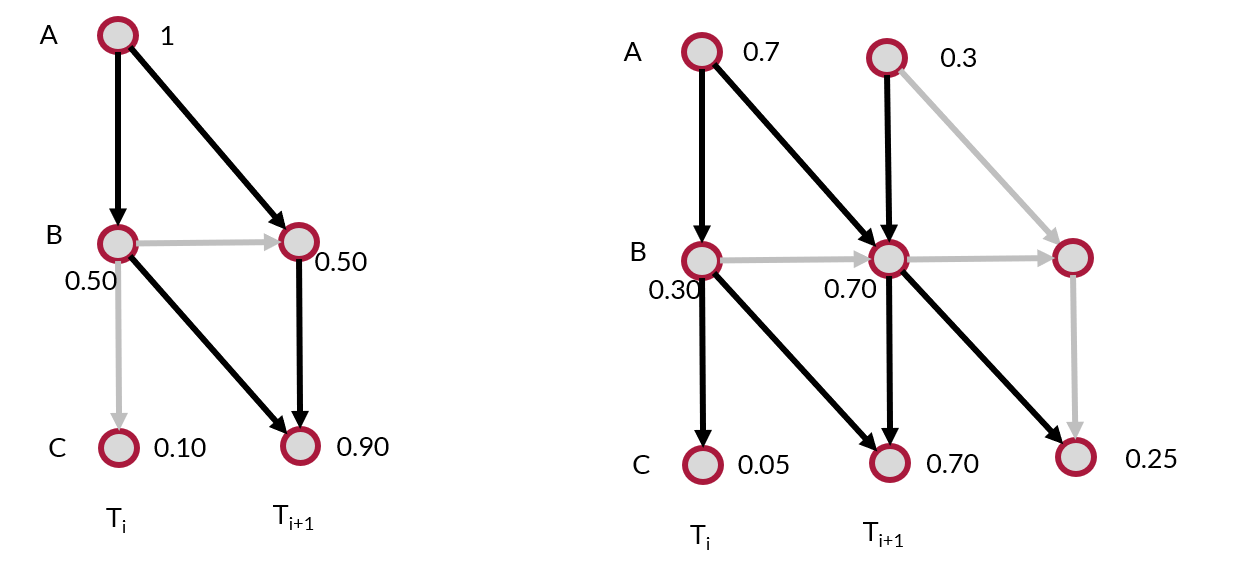}
\par\end{centering}
\caption{\label{fig:Flow-model-and}Flow model and volumes on nodes}
\end{figure}

Volumes that depart in the second time period must be known to not
be distributed over the whole second time period, shown with the blue
diagonal border lines for the volumes' travel along line A-B. The
corresponding flow is shown in figure \ref{fig:Flow-model-and}.

The aggregation constraints have to be reformulated in order to handle
volumes whose time windows stretches over several time periods. Note
that this could also be used to narrow the departure time windows
and/or the departure times within the time period. care must then
be taken, since capacity is only measured as a total amount for the
time period as a whole, and the flow model is not handling different
capacity usage within the time period. It it therefore essential to
choose the time period length in such a way that the model gives a
reasonable capacity resolution.

\section{Conclusions and future work}

We have presented a first description of a railway traffic flow model
based on a multi-commodity flow model. the model is innovative in
that the model handles volumes of trains rather than scheduling individual
trains. The model is also innovative in that time is discretized into
time periods and these periods together with the railway network forms
the flow model of arcs and nodes. On top of the traffic flow model
we introduce two layers, one to handle routes through the network
corresponding to Train Service Classes, TSCs \cite{Aronsson1302809}.
These TSCs package the possible services that can be run in the network.
A TSC is a named service from origin to the destination together with
important properties such as via stations, axle load, gauge, speed
category etc. The TSCs also contain a number (on or more) of named
routes through the network. The demand is distributed onto the named
routes according to congestion on the route's line segments and how
the objective function is formulated. The idea behind the TSCs is
that almost all of the demanded regular traffic on the railway network
should fit into declared TSCs. The third layer consist of the demands
for traffic, and should be matched on the TSCs to find the best match.
this matching also gives the departure time and transport duration
which in turn gives the earliest arrival at the destination. Different
routes can implement the same demand, thereby giving the route as
a decision variable in the model.

All together we have found in these first round of experiments that
the model is quite elaborate in terms of understanding the three layered
structure but computationally it delivers answers fast. We have tried
it with tests up to approximately 20 \% of the Swedish railway network
and are still able to deliver optimal answers within a couple of minutes.
Although these tests are encouraging there are more work to be performed
before we could say that the proposed model is ready for use. 

The objective function must be developed more carefully. The one presented
in this paper is too simple and do not have support in e.g. socio-economic
literature or other researched properties. Initial tests with an objective
function based on the Priority Criteria formulation of the Swedish
Network statement \cite{Trafikverket_JNB2021} has been performed.
It shows promising results and could be an important tool for constructing
and evaluating the capacity model, the capacity supply and as a traffic
valuation tool when investigating the impact of TCRs.

There are a number improvements to be addressed for the model to be
usable in practice. One already mentioned is to formulate the restriction
of postponed volumes, to get the demanded train volumes, required
to be integer, to depart as an evenly spread volume over time. Of
crucial importance is how the capacity properties are formulated so
that later timetabling can find a feasible schedule of train individuals
along with the proposed volumes from this model. The timetabling complexity
is ``hidden'' in the capacity formulation in the proposed model,
in terms of volumes of trains (or capacity units per train type allocated
to volumes). We have in this paper sketched how setup times could
be introduced in the model, as a function of the heterogeneity of
the traffic. Many industry sectors reason in such terms about future
production in the early planning phases, i.e. in terms of the number
of setups/changes on production facilities when producing different
products: how many setups that are needed, the duration of the setup
times between different product categories, the impact of longer series
of products on intermediate stores, warehouses etc. It is an interesting
question in itself if this can be done in railway traffic scheduling,
to form feasible high level plans with volumes of traffic usable in
e.g. calculating traffic supply and TCR impact in the early phases
of the capacity allocation process. The result from these high level
plans must have feasible timetables ``included'' in them, i.e. they
must be realizable if this approach shall be used to form a supply
of service classes to choose from. 

\printbibliography

\end{document}